\documentclass{article}      
\usepackage{amssymb}
\usepackage{amsmath}
\usepackage{latexsym}
\usepackage[mathscr]{euscript}
\usepackage{graphicx}
\usepackage{amscd}
\usepackage{amsthm} 
\usepackage{fullpage} 
\usepackage{wrapfig}
\newtheorem{theorem}{Theorem}
\newtheorem{corollary}[theorem]{Corollary}
\newtheorem{lemma}[theorem]{Lemma}  
\newtheorem{proposition}[theorem]{Proposition}

 \newtheorem{glue}[theorem]{Gluing Assumption}
\newtheorem{definition}{Definition}
 
\newcommand{\bc}{\mathbb{C}}

\newcommand{\bz}{\mathbb{Z}}
\newcommand{\br}{\mathbb{R}}

\newcommand{\p}{\partial}
\newcommand{\cc}{\mathcal{C}}
\newcommand{\ce}{\mathcal{E}}
\newcommand{\cj}{\mathcal{J}}

\newcommand{\cl}{\mathcal{L}}
\newcommand{\cs}{\mathcal{S}}

\newcommand{\cb}{\mathcal{B}}
\newcommand{\cp}{\mathcal{P}}

\newcommand{\cm}{\mathcal{M}}

\newcommand{\ct}{\mathcal{T}}

\newcommand{\ca}{\mathcal{A}}

\newcommand{\rz}{\br/\bz}
\newcommand{\hk}{\hookrightarrow}

\newcommand{\med}{\medskip}
\newcommand{\la}{\longrightarrow}
\newcommand{\bfl}{\begin{flushleft}}
\newcommand{\efl}{\end{flushleft}}

\newcommand{\eps}{\epsilon}

\newcommand{\calj}{\mathcal{J}}

\newcommand{\xr}{\xrightarrow}

 \newcommand{\bcm}{\bar \cm}

 \newcommand{\cah}{\ca_H}
 \newcommand{\ut}{\underbar{2}}
 \newcommand{\lkr}{\langle k \rangle}
 \newcommand{\ck}{\mathcal{K}}
 \newcommand{\trho}{\tilde \rho}

\begin{document}  

  \title{The Floer homotopy type of the cotangent bundle}
  \author{Ralph L. Cohen \thanks{The author was partially supported by a   grant  from the NSF} \\ Department of Mathematics  \\Stanford University \\ Stanford, CA 94305 }
 \date{\today}
\maketitle  
 \begin{abstract}  Let $M$ be a closed, oriented, $n$-dimensional manifold. In this paper we describe a spectrum in the sense of homotopy theory, $Z(T^*M)$,  whose homology is naturally isomorphic   to the Floer homology of the cotangent bundle, $T^*M$.    This Floer homology   is   taken with respect to a Hamiltonian
 $H : S^1 \times T^*M \to \br$ which is quadratic near infinity.   $Z(T^*M)$ is constructed assuming a basic smooth gluing result of $J$-holomorphic cylinders.  This spectrum will     have a $C.W$ decomposition with one cell for every  periodic solution of the  equation defined by the Hamiltonian vector field $X_H$.  Its induced cellular chain complex is exactly the Floer complex.   The attaching maps in the $C.W$ structure of 
 $Z(T^*M)$  are described in terms of the framed cobordism types of the moduli spaces
 of $J$-holomorphic cylinders in $T^*M$ with given boundary conditions.  This is done via a Pontrjagin-Thom construction, and an   important ingredient in this is proving, modulo this gluing result,  that these moduli spaces are compact, smooth, framed manifolds with corners.   We then prove that  $Z(T^*M)$, which we refer to as the ``Floer homotopy type" of $T^*M$, has the same homotopy type as the suspension spectrum of the free loop space, $LM$.  This generalizes the theorem first proved by C. Viterbo that  the Floer homology of $T^*M$ is isomorphic to   $H_*(LM)$. 
     \end{abstract}

 \tableofcontents

 \section*{Introduction} 
 An intriguing theorem proved by Viterbo \cite{viterbo}, with alternative proofs by Salamon and Weber \cite{salamonweber},  as well as Abbondondolo and Schwarz \cite{abboschwarz},  states that the Floer homology of the cotangent bundle of a closed, oriented manifold $M$,  is isomorphic to the homology  of the free loop space, 
$$
 HF_*(T^*M) \cong H_*(LM).
$$
 Here $T^*M$ is endowed with its canonical symplectic structure, and the Floer homology is taken with respect to a time dependent Hamiltonian,
$$
 H : \rz \times T^*M \to \br
$$
 which is quadratic near infinity in an appropriate sense.  The goal of this paper is to 
 examine the homotopy theoretic underpinnings of this isomorphism.
 
 Recall that the Floer homology of any symplectic manifold, $(N, \omega)$ is defined to be the homology of the Floer complex, $CF_*(N)$ which is generated by the $1$-periodic solutions  of a Hamiltonian equation of the form
\begin{equation}\label{ham}
 \frac{dx}{dt} =  X_H(t, x(t)),
\end{equation}
where  $H$ is a time dependent Hamiltonian, and  $X_H$ is the   Hamiltonian vector field.   

Floer theory can be viewed as a generalized Morse theory, but in the case of a classical Morse  function $f : P \to \br$, where $P$ is a smooth, closed manifold,  there is more information than just the corresponding chain complex.  Namely, the  Morse complex    $C^f_*(P)$ can be viewed as  the cellular chain complex of a  $C.W$ complex $Z_f(P)$ which is naturally homotopy equivalent to the manifold $P$, and has one cell for each critical point of the function $f$.  

It is then natural to  ask, in the case of Floer theory, is there a naturally defined underlying  $C.W$ complex (or spectrum)   with one cell for each solution
to the Hamiltonian equation (\ref{ham}), so that its associated cellular chain complex
is isomorphic to the Floer complex, $CF_*(N)$?  This question was taken up by the author, Jones, and Segal in \cite{cjs} where conditions on a ``Floer functional" were obtained that allowed the construction of such an underlying homotopy type.  The most important of these conditions, which significantly restricts when such a Floer homotopy type can exist, is that the moduli space of gradient flow lines connecting
two critical points, have a compactification which is a framed manifold with corners.  Moreover the framings must be chosen compatibly with respect to the gluing of 
these moduli spaces. 

The goal of this paper is to show that, modulo a specific analytic gluing construction,  these conditions are satisfied by the symplectic action functional on the loop space of $T^*M$,  perturbed by an appropriate Hamiltonian, $H : \rz \times T^*M  \to \br$,   
$$
\ca_H : L(T^*M) \to \br.
$$
We then explicitly describe the resulting Floer homotopy type of $T^*M$, which we denote
by $Z(T^*M)$.  In this case $Z(T^*M)$ is a $C.W$-spectrum with one cell for every periodic solution
of the Hamiltonian equation (\ref{ham}).   

To be more specific, let $J$ be a compatible almost complex structure on $T^*M$, let $a$, $b : S^1\to T^*M$ be $1$-periodic solutions of the Hamiltonian equation, and let $\bcm (a,b; H,J)$ be the space of piecewise $J$-holomorphic cylinders $u: S^1 \times \br \to  T^*M$ which converge to $a$ as $t \to -\infty$, and $b$ as $t\to +\infty$.  (The precise definitions of these moduli spaces will be given below.)  The smoothness,  compactness, and orientability  properties of the uncompactified moduli spaces $\cm(a,b; H,J)$ were studied by Abbondandolo and Schwarz in \cite{abboschwarz}. The topology of the compactified moduli space will be discussed in section 2.  The basic gluing result that we will assume in this paper is the following.

 \begin{glue}  \label{glue}    Let $a_0, \cdots, a_n$ be a sequence of $1$-periodic solutions so that $\cm (a_i, a_{i+1}; H,J)$ is nonempty  or each $i$.  Then there is an $\eps >0$ and   local diffeomorphisms (i.e a diffeomorphisms onto   open subsets)
$$
\cm(a_0, a_1) \times \cdots \cm (a_{n-1}, a_n) \times [0, \eps)^{n-1} \to \bcm (a_0, a_n)
$$ 
which give $\bcm(a_0, a_n)$ the structure of a smooth, compact manifold with corners.  In particular,  the boundary $\p\bcm(a_0, a_n)$ consists of the images of these gluing maps restricted to $\cm(a_0, a_1) \times \cdots \cm (a_{n-1}, a_n) \times \{0\}$, and on the open intervals, these gluing maps restrict to give local diffeomorphisms,
$$
\cm(a_0, a_1) \times \cdots \cm (a_{n-1}, a_n) \times (0, \eps)^{n-1} \to \cm (a_0, a_n).
$$  
\end{glue}
 
\med
\noindent
\bf Remark.  \rm  Gluing constructions of the above type have been constructed in the several places in the literature (see for example \cite{floer}, \cite{taubes}, \cite{donaldson},  \cite{baraudcornea}).   However  specific gluing results that would imply  
that the compact moduli spaces $\bcm (a,b; J,H)$ studied here are smooth manifolds with corners, have not appeared in the literature. The analysis required to prove such a gluing result  is of a very different sort of mathematics than the algebraic topological techniques that are used in this paper, and will not be pursued here.    Thus the statements in this paper can be viewed as topological results that would follow from this analysis.  Although we state these results as theorems, it should be understood that for their proofs to be complete, the   above gluing assumption must be proved. This, as far as the author understands,   has not yet been completed.  From now on in this paper we will operate under  Gluing Assumption \ref{glue}.

 Our first  result states that  assuming the above gluing result,  these moduli spaces have natural framings. 

\begin{theorem}\label{frameable}  For appropriate choices of Hamiltonian $H$, and with respect to a generic choice of almost complex structure $J$ on $T^*M$, the moduli spaces $\bcm (a,b; H,J)$ are smooth, compact, framed manifolds with corners. The dimension of this moduli spaces is $\mu(a) - \mu (b) -1$, where $\mu$ represents the Conley-Zehnder index. 
\end{theorem}

\med
We will actually prove something stronger than this.  If $\cc_H$ is the ``flow category", of $\ca_H : L(T^*M) \to \br$, whose objects are the critical points (i.e periodic  solutions to the Hamiltonian equation), and whose space of morphisms  between two solutions $a$ and $b$ is the moduli space $\bcm(a,b; H, J)$, then when certain conditions are met,   $\cc_H$   is a ``compact, smooth, framed category of Morse-Smale type".   These conditions will be defined below, but basically they refer to a certain transversality conditions, and the fact that  the choices of framings are compatible with the composition in the category, which in this case corresponds to gluing of flow lines. 

These conditions  were what was required in   \cite{cjs}
to produce a ``Floer homotopy type".  In the general setting studied in \cite{cjs}, the Floer homotopy type was a certain inverse limit of spectra, known as a ``prospectrum".  In the case of the cotangent bundle,   however,  since    the critical values of $\ca_H$ are bounded below,  we will observe   that  the   procedure in \cite{cjs} actually produces a $C.W$-spectrum $Z(T^*M)$ which realizes the Floer homotopy type.  Namely, we prove the following.

\begin{theorem}\label{floerhomotopy} For appropriate choices of Hamiltonian $H$, and with respect to a generic choice of almost complex structure $J$ on $T^*M$, there is an associated finite type $C.W$-spectrum $Z(T^*M)$ with one cell for every critical point of the perturbed sympectic action functional, $\ca_H : L^(T^*M) \to \br$, i.e solutions to the Hamiltonian equation, $ \frac{dx}{dt} =  X_H(t, x(t))$.  The attaching maps in this cellular structure are defined explicitly  using the Pontrjagin-Thom construction, in terms of the framed bordism classes of manifolds with corners represented by  the moduli spaces, $\bcm(a,b; H, J)$. The induced cellular chain complex of $Z(T^*M)$ is precisely the Floer complex taken with respect to this Hamiltonian and almost complex structure,  $CF_*(T^*M; H, J)$. 
\end{theorem}

\med
This theorem says that the Floer homology, $HF_*(T^*M)$ is computed by taking the homology of the complex obtained by applying ordinary homology to the subquotients of the cellular filtration of $Z(T^*M)$.   These subquotients are wedges of spheres indexed by the critical points of $\ca_H$, and his theorem says that this complex is the Floer complex, $CF_*(T^*M)$.    Notice that if one applies a \sl generalized \rm homology theory $h_*$ to these subquotients, one obtains a spectral sequence converging to the ``Floer $h_*$-theory" of $T^*M$.

\begin{corollary}\label{generalized}.  Let $h_*$ be any generalized homology theory.
There is a a spectral sequence whose $E_1$-term is given by
$$
E_1^{p,q} = CF_p(T^*M) \otimes h_q(point)
$$
and converges to $h_*(Z(T^*M))$, the ``Floer $h_*$-theory of $T^*M$". Here $p$ is the Conley-Zehnder index of the periodic solution to the Hamiltonian.    In the case when $h_* = H_*$, ordinary homology, then the $E_1$ -term is exactly the Floer complex,
and the spectral sequence collapses at the $E_2$-level, which is the Floer homology, $HF_*(T^*M)$. 
\end{corollary}

\med
In this spectral sequence, the differential at the $E_1$ level, $$d_1 : CF_p(T^*M) \otimes h_q(point) \to CF_{p-1}(T^*M) \otimes h_q(point)$$ is given  in terms of the numbers $\#\bcm(a,b; H,J)$ when $a$ has Conley-Zehnder index $\mu (a) = p$, and $\mu (b) = p-1.$ The number, $\#\bcm(a,b; H,J)$, which is given by counting the number of points in $\bcm(a,b; H,J)$, with signs determined by orientations, can also be viewed as the framed bordism class of this zero dimensional   manifold. (In dimension zero a framing and an orientation are the same thing.)  The higher differentials in this spectral sequence are determined by the framed bordism classes of the higher dimensional moduli spaces. In the case of ordinary homology, the collapse of this spectral sequence implies
that one does not need to consider these higher dimensional spaces.  However for a generalized homology theory  these higher dimensional moduli  spaces   play a critical role.  

\med
Our final result is an identification of the Floer homotopy type, $Z(T^*M)$, of the cotangent bundle.   We compare the underlying homotopy theory of a Morse function on the loop space, $\ce : LM \to \br$, with the Floer homotopy type $Z(T^*M)$.  This involves comparing the framed bordism types of the moduli spaces of gradient trajectories of the Morse function $\ce$, with those of the moduli spaces of $J$-holomorphic cylinders $\bcm(a,b; H, J) $.   This adapts methods of  Abbondandolo and Schwarz \cite{abboschwarz},  to the setting of framed manifolds with corners. The result of this study is the following.

  Given a space $X$,  let $\Sigma^\infty (X_+)$ denote the suspension spectrum of $X$ with a disjoint basepoint.  

\med
\begin{theorem}\label{equivalence} Given the hypotheses of Theorem \ref{floerhomotopy}, there is a homotopy equivalence of spectra,
$$
\Phi : \Sigma^\infty (LM_+) \xr{\simeq}  Z(T^*M).
$$
 \end{theorem}
 
 Notice that this  generalizes the Viterbo theorem stating that $HF_*(T^*M) \cong H_*(LM)$.  
 
 \med
The motivation for this work came from recent work by several authors whose goal is to relate the string topology theory of $LM$, as originally defined by Chas and Sullivan \cite{chassull}, to  the symplectic   theory of $T^*M$ (see for example  \cite{abboschwarz2}, \cite{fukaya}, \cite{cohenmont}, \cite{cielelai}).  These approaches require the string topology operations to be defined and have good properties on the chain level.  Our goal is to 
take a more homotopy theoretic approach.  In \cite{cj}, and \cite{cg} the author, Jones, and Godin showed that string topology operations can be realized on the homotopy level.  The results of this paper, in particular Theorems \ref{floerhomotopy} and \ref{equivalence},  realize the Floer theory of $T^*M$ on the homotopy level.  This program was announced and summarized in \cite{cohenmont}.  

\med
This paper is organized as follows.  In section one we describe the basic homotopy theory that is necessary.  In particular  we describe  conditions that completely characterize the realizations of   a finite chain complex of free abelian groups,   as the cellular chains of a finite $C.W$ complex or spectrum.  These conditions are also described in terms of cobordism classes of framed manifolds with corners.  This discussion is  a recasting of the discussion in \cite{cjs}, and may be of independent interest.  (See Theorem \ref{realize} below.)
In section two we study the Floer theory of $T^*M$, and prove that the corresponding flow category $\cc_H$ satisfies the conditions necessary for generating   a stable homotopy type (``Floer homotopy type") as described in section one. We    then prove Theorems \ref{frameable} and \ref{floerhomotopy}. In section three we  identify the Floer homotopy type with the stable homotopy type of the free loop space  thereby proving Theorem  \ref{equivalence}.   

\med
The author is very grateful to O. Cornea, Robert Lipshitz, C. Manolescu, and K. Wehrheim  for very helpful discussions and correspondence about this material, and in particular for helping the author understand the status of gluing constructions in symplectic geometry.  

\section{The homotopy theory}

From the algebraic topology point of view, the question of  naturally realizing the Floer chain complex  by an underlying homotopy type, is a special case of the question of  understanding how a given  chain complex
$$
\to \cdots \to C_i \xr{\p_i} C_{i-1} \xr{\p_{i-1}} \cdots \to C_0
$$
 may be realized as  the cellular chain complex  of a $C.W$-complex or spectrum.  This question was addressed   in \cite{cjs}.  In this section we  recall that discussion and give functorial criteria.  We then use Pontrjagin-Thom theory
to recast these criteria in terms of framed manifolds with corners.

 \subsection{Stable attaching maps of finitely filtered spaces }
 
 Recall in \cite{cjs} the authors described how, given a compact space $X$,  filtered by compact subspaces,
$$
X_0 \hk X_1 \hk \cdots \hk X_n = X,
$$
where each $X_{i-1} \hk X_i$ is a cofibration with cofiber, $K_i = X_i/X_{i-1}$,
then one can ``rebuild" the homotopy type of the $n$-fold suspension, $\Sigma^n X,$
as the union of iterated cones and suspensions of the $K_i$'s,
\begin{equation}\label{decomp}
\Sigma^n X \simeq \Sigma^nK_0 \cup c(\Sigma^{n-1} K_1) \cup \cdots \cup c^i(\Sigma^{n-i} K_i) \cup \cdots \cup c^n K_n.
\end{equation}

This decomposition can be described as follows.   The cofibration sequences $X_{i-1} \xr{u_{i-1}} X_{i } \xr{p_{i }} K_{i }$ have Puppe extensions,
$$
K_i \xr{v_i} \Sigma X_{i-1} \xr{u_{i-1}} \Sigma X_i.
$$
Let $\p_i : K_i \to \Sigma K_{i-1}$ be the composition $\p_i = \Sigma p_{i-1}\circ  v_i  : K_i \to \Sigma X_{i-1} \to \Sigma K_{i-1}$.    Consider the following sequence:

\begin{equation}\label{htpychain}
K_n \xr{\p_n} \Sigma K_{n-1} \xr{\p_{n-1}} \cdots \xr{\p_{i+1}}\Sigma^{n-i }K_{i } \xr{\p_{i}} \Sigma^{n-i+1}K_{i-1} \xr{\p_{i-1}} \cdots \xr{\p_1}\Sigma^{n}K_0 = \Sigma^n X_0.
\end{equation}
In this sequence, for ease of notation, $\p_j$ represents the appropriately iterated suspension of the map $\p_j$ defined above. 
This sequence  can be viewed as a ``homotopy chain complex"  because each of the compositions,
$$
\p_j \circ \p_{j+1}
$$
is canonically null homotopic.  This is because the composition contains as a factor, the cofibration sequence,
$\Sigma^{n-j}X_j \xr{\Sigma^{n-j}p_j} \Sigma^{n-j}K_j \xr{\Sigma^{n-j}v_j} \Sigma^{n-j+1}X_{j-1}$.   This canonical null homotopy
defines an extension of $\p_j : \Sigma^{n-j}K_j \to \Sigma^{n-j+1}K_{j-1}$ to the mapping cone,
$$
c(\Sigma^{n-j-1}K_{j+1}) \cup_{\p_{j+1}} \Sigma^{n-j}K_j   \la \Sigma^{n-j+1}K_{j-1}.
$$
More generally for every $q$, using these null homotopies,  $\p_j : \Sigma^{n-j}K_j \to \Sigma^{n-j+1}K_{j-1}$  extends to an attaching map for   map of the iterated mapping cone,
\begin{equation}\label{attach}
c^q(\Sigma^{n-j-q}K_{j+q}) \cup c^{q-1}(\Sigma^{n-j-q+1}K_{j+q-1}) \cup \cdots \cup c(\Sigma^{n-j-1}K_{j+1}) \cup_{\p_{j+1}} \Sigma^{n-j}K_j  \la  \Sigma^{n-j+1}K_{j-1}. 
\end{equation}

To keep track of the combinatorics of these attaching maps,  a category $\calj$ was introduced in \cite{cjs}.   The objects of $\calj$ are the integers, $\bz$.   To describe the morphisms, we first define,   for any two integers $n>m$,  the space 
\begin{equation}
J(n,m) = \{t_i, \, i \in \bz, \quad \text{where each $t_i$ is a nonnegative real number, and}   \quad t_i=0,  
 \, \text{unless} \, m < i < n. \}
\end{equation}
Notice that $J(n,m) \cong  \br_+^{n-m-1}$, where $\br_+^q $ is the space of $q$-tuples of   nonnegative real numbers.     Notice that one has natural inclusions,
$$
\iota: J(n,m) \times J(m,p)  \hk J(n,p).
$$
\med
We then define the  morphisms   in $\calj$  as follows.  For integers $n<m$ there are no morphisms from $n$ to $m$.  The only morphism from an integer $n$ to itself is the identity. If $n = m+1$, we define the morphisms to be the two point space, $ Mor(m+1, m) = S^0$.   If       $n>m+1$,  $Mor(n,m)$  is given by the one point compactification,
$$
Mor (n,m) =  J(n,m)^+ = J(n,m) \cup \infty.
$$
For consistency of notation we refer to all the morphism spaces  $Mor(n,m)$ as $J(n,m)^+$.  
Composition in the category is given by addition of sequences,
$$
J(n,m)^+ \times J(m,p)^+ \to J(n,p)^+.
$$ 
If $n-m \geq 2$, notice that   $Mor(n,m) = J(n,m)^+$ is topologically a   disk of dimension $n-m-1$ with a distinguished basepoint $(=\infty)$. Notice also that for a based space $Y$,  the smash product $J(n,m)^+\wedge Y$  is the iterated cone,
$$
J(n,m)^+\wedge Y = c^{n-m-1}(Y).
$$

Given integers $p > q$, then there are subcategories  $\calj^p_q$ defined to be  the full subcategory generated by integers $q \geq m \geq p.$
Notice that given a filtered space $X = X_n$ as above, there is an induced functor,
$$
Z_X : \calj^n_0 \to  Spaces_*
$$
where $ Spaces_*$ is the category of compact spaces with basepoint.  On objects we have
$$
Z_X(m) = \Sigma^{n-m}K_m  =  \Sigma^{n-m}(X_m/X_{m-1}).
$$  
On morphisms, $Z_X$ is defined via the relative attaching maps,
$$
Z_X : J(m,p)^+ \wedge \Sigma^{n-m} K_m = c^{m-p-1}( \Sigma^{n-m} K_m) \la  \Sigma^{n-p}K_p
$$
 given in  (\ref{attach}) above.    A precise description of this functor was given in \cite{cjs}. 
 
      As described in \cite{cjs},  given such a functor,
 $  Z : \calj^p_q \to  Spaces_*$  one can take its geometric realization,
\begin{equation}\label{zrealize}
 |Z| = \coprod_{q\leq j \leq p} Z(j)\wedge J(j, q-1)^+ / \sim
\end{equation}
 where one identifies the image of $Z(j) \times J(j,i)^+\times J(i, q-1)^+$  in $Z(j) \wedge J(j, q-1)^+$  with its image in $Z(i) \wedge J(i, q-1)$ under the map on morphisms.  
 Notice that since  $Z(j) \wedge J(j, q-1)^+ $ is the iterated cone  $ c^{j-q}Z(j)$,  
\begin{equation}\label{cell}
 |Z| = Z(q) \cup c(Z(q+1)) \cup 
 \cdots \cup c^{p-q}(Z(p)).
 \end{equation}
 There is a double complex, 
\begin{equation}\label{complex}
 C_*(Z) = \oplus_{q \leq j \leq p} C_*(Z(j))
\end{equation}
 which computes the homology of $|Z|$.  
 In the case of the functor $Z_X: \calj^n_0 \to  Spaces_*$ coming from a filtered space as above, this decomposition is exactly the decomposition given in (\ref{decomp}), and
 $|Z_X| \simeq \Sigma^n X$.   
 
 Notice that in this construction, our functor $Z_*$ might just as well have taken values in the category $Sp_*$ of finite spectra. Here the identifications in the definition (\ref{zrealize}) of $|Z|$, would be replaced by coequalizers in $Sp_*$ in the usual way. 
 
 We observe that the above argument proves the following theorem which gives our basic criterion for realizing a finite chain complex as a homotopy type.
 
 \med
 \begin{theorem}\label{realize}  Each realization of a  finite chain complex of finitely generated free abelian groups,
 $$
 C_n \xr{\p_n} C_{n-1} \to \cdots \to C_i \xr{\p_i} C_{i-1} \xr{\p_{i-1}} \cdots \xr{\p_1} C_0
 $$
 as the  the cellular chain complex of a finite $C.W$ spectrum $X$, with
 $C_i = H_i(X^{(i)}, X^{(i-1)})$, occurs as the geometric realization of a functor  
   $$
 Z_X :  \calj^n_0 \to  Sp_*
 $$
 with $Z_X(i) = \Sigma^{n-i}(X^{(i)}/ X^{(i-1)})$.
 \end{theorem}

 \subsection{Framed manifolds with corners}
Consider  the case when  a functor $Z : \calj \to  Sp_*$ has the property that  each $Z(i)$ is homotopy equivalent to a wedge of spheres (i.e a wedge of suspension spectra of spheres).    This, of course, is the case when $Z$ is induced by a cellular filtration of a space.  In this setting, the   maps on morphisms  are defined   by   collections of maps  
 $$
 Z(\alpha, \beta ) : S^k_\alpha \wedge J(i,j)^+ \la S^r_\beta
 $$ for each sphere $S^k_\alpha$ in the wedge decomposition of $Z(i)$ and $S^r_\beta$ in the wedge decomposition of $Z(j)$.  We can find a smooth representative of this map on the level of spaces, $ Z(\alpha, \beta ) : S^{k+L}_\alpha \wedge J(i,j)^+ \la S^{r+L}_\beta$.  Then  the inverse image of a regular point   is a compact, smooth manifold with corners, $M$, embedded
 in $\br^{k+L} \times J(i,j) \cong \br^{k+L} \times \br_+^{i-j-1}$, where the embedding respects the corner structure and has a trivialized normal bundle.  By Pontrjagin-Thom theory, the ``framed cobordism type" of this manifold is determined by and determines this attaching map.  The purpose of this section is to describe categorically,  the condition on a collection of framed manifolds
 with corners, so that it determines a functor $Z : \calj \to  Sp_*$, and therefore a stable homotopy type.   
  
 \med
 Recall that an $n$-dimensional manifold with corners $M$  has charts which are local homeomorphisms with  $\br_+^n$. 
 Recall that the boundary  $\p M$ are those points in $M$ which map under a local chart to the boundary   $\p \br_+^n = \{(t_1, \cdots , t_n) \in  \br_+^n \,\, \text{such that at least one of the $t_i$'s is zero.}\}$.     Let $\psi : U \to (\br_+)^n$ be a chart of a manifold with corners $M$.  For $x \in U$,  the number of zeros of this chart, $c(x)$ is independent of the chart.     A connected face, or ``boundary hypersurface" of a manifold with corners $M$ is the closure of a component of $\{x\in M \, : \, c(x) = 1\}$.  
  Recall from \cite{janich}, \cite{laures}  that a \sl manifold with faces \rm is a manifold with corners $M$  such that each $x \in M$  lies in $c(x)$ different, connected faces.  Notice that in a manifold with faces, any disjoint union of connected faces is itself a manifold with faces.
  
  Recall furthermore, that a  $\langle k \rangle$-manifold  $M$, is a manifold with faces, together with an ordered $k$-tuple $F_1(M), \cdots, F_k (M)$ of faces, satisfying
  \begin{enumerate}
  \item The union of these faces is the entire boundary,
  $$F_1(M) \cup \cdots \cup F_k(M) = \p M.$$
  \item  Each intersection   $F_i (M) \cap F_j(M)$ is a  face of both $F_i(M)$ and of $F_j(M)$.
   \end{enumerate}
   
 The  archetypical example of a  $\langle k \rangle$-manifold is $ \br_+^k$.  In this case
 the face $F_j \subset \br_+^k$ consists of those $k$-tuples with the $j^{th}$- coordinate equal to zero.
Cobordisms of   $\langle k \rangle$-manifolds have been studied by Laures in \cite{laures}, as have their basic embedding properties.   We will make considerable use of these properties in this paper.   

For example as Laures indicated, the data of a $\langle k \rangle$-manifold can be encoded in a   categorical way  as follows.  Let $\underbar{2}$ be the partially ordered set with two objects, $\{0, 1\}$, generated by a single nonidentity morphism $0 \to 1$.  Let $\ut^k$ be the product of $k$-copies of the category $\ut$.  A $\langle k \rangle$-manfold $M$ then defines a functor from $\ut^k$   to the category of topological spaces,  where for an object  $a = (a_1, \cdots , a_k) \in \ut^k$,
$M(a)$ is the intersection of the faces $F_i(M)$   with $a_i = 0$.  Such a functor is a $k$-dimensional cubical diagram of spaces, which, following Laures' terminology, we refer to as a $\langle k \rangle$-diagram.  Notice that $\br_+^k(a) \subset \br_+^k$ consists of those $k$-tuples of nonnegative real numbers so that the $i^{th}$-coordinate is zero for every $i$ such that $a_i=0$.

    In this section we will consider embeddings of manifolds in corners into Euclidean spaces of the form given by the following definition.

\begin{definition}\label{embed} An embedding of a $\langle k \rangle$-manifold $M$ into Euclidean space $\br^m \times  \br_+^k$  is a natural transformation of $\langle k \rangle$-diagrams
$$
e : M \hk \br^m \times  \br_+^k
$$ for some $m$, that satisfies the following properties:
\begin{enumerate}
\item For every object $a \in \ut^k$, $e(a) : M(a) \to \br^m \times  \br_+^k(a)$ is an inclusion of a smooth submanifold, and 
\item for all $a > b$,  the intersection  $M(a) \cap  \left( \br^m \times \br_+^k(b)\right)$ is equal to $M(b)$. 
\end{enumerate}
 \end{definition}

\med
Such an embedding was called a \sl neat \rm embedding by Laures in \cite{laures}.  Moreover he proved that  a manifold with corners $M$ admits such an embedding if and only if it is a  $\langle k \rangle$-manifold.  In particular if $e : M \hk \br^m \times  \br_+^k$  is such an embedding,  the $\langle k \rangle$-structure on $M$ is given by the intersection,  $F_j (M) = M \cap \left( \br^m \times F_j( \br_+^k)\right).$ Because of this, we refer to such an embedding of $\langle k \rangle$-manifolds simply as an embedding of manifolds with corners. 

\med
Notice that given  an embedding of manifolds with corners,  $e : M \hk \br^m \times  \br_+^k$, then it has a well defined normal bundle.  In particular, for any pair of objects in $\ut^k$ $a > b$, then the normal bundle of $e(a) : M(a) \hk \br^m \times \br_+^k(a)$, when restricted to $M(b)$, is the normal bundle of $e(b) : M(b) \hk \br^m \times \br_+^k(b)$. 

 Said another way, the normal bundle is classified by a homotopy class of maps  (natural transformations) of  $\langle k \rangle$-diagrams
 $\nu_e : M \to BO(q)$,  where $BO(q)$ represents the constant $\lkr$-diagram whose value at every object is a model of the classifying space $BO(q)$,  and whose arrows are all equal to the identity map. (By a ``model" we mean a choice of space of the given homotopy type.)  Here $q$ is the codimension, $n = m+k- \,dim \, M$. By a homotopy   of maps of $\lkr$-diagrams we mean
 a continuous, one-parameter family of such maps, in the usual way. 
 
 Similarly, the \sl stable normal bundle \rm of $M$ is classified by a homotopy class of maps of $\lkr$-diagrams
 $$
 \nu_M : M \to BO
 $$
 for some model of $BO$, viewed as a constant $\lkr$-diagram.  This homotopy class is independent of any choice of embedding into Euclidean space.

For our purposes, we would like to consider ``framed" embeddings of manifolds with corners, which has the extra structure of a trivialization of the normal bundle.

\begin{definition}\label{framedcorner}  Let $M$ be a $\lkr$-manifold.  Let $\nu_M : M \to BO$ be a classifying map of its stable normal bundle. That is, $BO$ is a constant $\lkr$-diagram whose value on an object is a model for the classifying space $BO = \lim_{n\to \infty}BO(n)$, and $\nu_M: M \to BO$ is a map of $\lkr$-spaces which on every object   classifies the stable normal bundle of the underlying manifold.  Let $p :EO \to BO$ be a fibration, where $EO$ is contractible.  View $EO$ as a constant $\lkr$-diagram, and $p : EO \to BO$ a map of $\lkr$-diagrams.  Then a \sl framing \rm of the stable normal bundle is a homotopy class of lifting 
$$
\tilde \nu_M : M \to EO
$$
of $\nu$.
\end{definition}

Now consider again   an embedding of manifolds with corners,   
$$
e : M \hk \br^m \times  \br_+^k.
$$ of codimension $q$.  
If  $q$ is sufficiently large and  $M$ has a framing of its stable normal bundle,    then a choice of  framing  defines, via the tubular neighborhood theorem,  an extension of $e$ to an embedding
\begin{equation}\label{framedemb}
e :  \br^q \times  M \to  \br^m \times (\br_+)^n
\end{equation}
 which is a diffeomorphism onto an open neighborhood of the image $e(M)$.  In particular on every object $a \in \ut^k$, this map   restricts to give a local diffeomorphism
$e :  \br^q \times M(a)   \to \br^m \times   \br_+^k (a) $.   This extension is well defined up to a self diffeomorphism of $\br^q \times M$, fixing $\{0\} \times M$,  which is isotopic to the identity.  We refer to this property as being   ``well defined up to \sl source isotopy".  \rm 
 We call  such a class of embedding
with extension as a ``framed"   embedding of a manifold  with corners. 

\med
Notice that given a framed embedding, $e : \br^q \times  M \to  \br^m \times  \br_+^k$,
the Pontrjagin-Thom construction gives us a map from the one point compactifications,
$$
\tau_e : S^m \wedge (\br_+^k \cup \infty) \to (\br^q \times M)\cup \infty \xr{project} \br^q \cup \infty = S^q.
$$
The source of this map is the iterated cone, $c^k (S^m)$, or equivalently the space
$S^m \wedge J(i,j)^+$ when $i = j+k+1$.  So the Pontrjagin-Thom  construction gives a map
\begin{align}\label{pont-thom}
\tau_e : \, &c^k(S^m)   \to S^q \\
&S^m \wedge J(i,j)^+  \to  S^q. \notag
\end{align} Here again, $m + k - q = dim \, M$.  
Conversely, given a smooth map, $\tau : c^k(S^m)   \to S^q $ then the inverse image
of a regular point is a compact, smooth $\lkr$-manifold $M$ of dimension $k+m-q$, equipped with a framed embedding,
$e : \br^q \times  M \hk  \br^m \times  \br_+^k$.  

\subsection{Compact, smooth, framed categories}
From the point of view of Theorem \ref{realize}  and the discussion in section 1.2 , in order
to define a stable homotopy type one needs an appropriately compatible collection of framed manifolds with corners.  In this section
we make this precise by defining the notion of a ``compact, smooth, framed category".
 This is a slight variant of the notion defined in \cite{cjs}.    
 
 \med
 We adopt the following definition.

\begin{definition}\label{framed}(\cite{cjs})A smooth, compact  category  is a topological category $\cc$ whose objects form a discrete set, and whose whose morphism spaces, $Mor (a,b)$ are compact, smooth manifolds with corners,  such that the composition maps $\mu_ : Mor (a,b) \times Mor (b, c) \to Mor(a, c)$ are smooth codimension one embeddings (of manifolds with corners) into the boundary.  

A smooth, compact category $\cc$  is said to be a ``Morse-Smale" category if the following
additional properties are satisfied.  
\begin{enumerate}
\item The objects of $\cc$ are partially ordered by condition
$$
 a \geq  b   \quad \text{if} \quad Mor (a,b) \neq \emptyset.
$$
\item $Mor(a,a) = \{identity \}$. 
\item  There is a set map, $\mu : Ob (\cc) \to \bz$ which preserves the partial ordering
so that  if $a > b$, 
$$
dim \, Mor(a,b) = \mu (a) - \mu (b) -1.
$$
$\mu$ is known as an ``index" map.   A Morse-Smale category such as this  is said to have finite type, if for each pair of objects  $a > b$,  there are only finitely many objects
$\alpha$ with $a > \alpha > b$. 
\end{enumerate}
\end{definition}

 We are now in a position to define a ``framing" of such a smooth category.
 
 \med
 
 \begin{definition}\label{framedembed}  Let $\cc$ be a compact, smooth, Morse-Smale category of finite type.
Let $a > b$ be objects, and let $\cc^a_b \subset \cc$ be the full subcategory generated by objects $\alpha$ with $a> \alpha > b$.   By a framed embedding of $\cc^a_b$  we mean that
for every pair of objects $\alpha, \beta$ in  $\cc^a_b$ with  $\alpha > \beta$,    there is a framed embedding of manifolds with corners,
 $$
 e_{\alpha, \beta} : Mor (\alpha, \beta) \subset  Mor (\alpha, \beta) \times \br^L\hk \br^L \times J(\mu (\alpha), \mu (\beta))
 $$
 satisfying the compatibility requirement that given any three objects, $\alpha > \beta > \gamma $ in $\cc^a_b$,  then the following diagram of embeddings  commutes up to source isotopy:
  $$
\begin{CD}
Mor(\alpha, \gamma) \times \br^L @>e_{\alpha, \gamma}>\hk >  \br^L \times J(\mu(\alpha), \mu (\gamma) \\
@AAA    @AA 1 \times \iota A  \\
@A compose A A        \br^L \times  J(\mu (\alpha), \mu (\beta))  \times  J(\mu(\beta), \mu (\gamma)) \\
@AAA      @A  A  e_{\alpha, \beta} \times 1 A \\ 
Mor (\alpha, \beta) \times Mor(\beta, \gamma) \times \br^L @>\hk>1 \times e_{\beta, \gamma} >  Mor (\alpha, \beta) \times \br^L \times J(\mu(\beta), \mu (\gamma)).
    \end{CD}
$$   In particular, this means that each morphism space $Mor (\alpha, \beta)$ is a $\langle \mu (\alpha)-\mu (\beta) -1\rangle$-manifold. 

A framing of a compact, smooth, Morse-Smale category of finite type is a choice of  framed embedding of $\cc^a_b$ for each pair of objects, $a > b$. In particular, this means that each morphism space $Mor (\alpha, \beta)$ is a $\langle \mu (\alpha)-\mu (\beta) -1\rangle$-manifold. 
\end{definition}

\med
The following was essentially proved in \cite{cjs}. (We say ``essentially" because some of the language of \cite{cjs} is slightly different than what is used here.)

\med
\begin{theorem}\label{fltype}  Let $\cc$ be a compact, smooth, framed Morse-Smale category of finite type.  For objects $a > b$,  Let $p = \mu (a)$, and $q = \mu (b)$ be the indices.   Then  there is a functor $Z^{a,b}_{\cc} : \calj^p_q \to  Sp_*$, whose geometric realizations, $|Z^{a,b}_{\cc}|$ fit together to give a prospectrum.  This prospectrum is called the   ``Floer homotopy type" of $\cc$.  
\end{theorem}

\begin{proof}  For completeness, we describe the functors $Z^{a,b}_{\cc}: \calj^p_q \to  Sp_*$.  We will not review the prospectrum aspect of this theory, since we will not need it   in this paper.   Let $\cc$ be such a compact, smooth, Morse-Smale category, and let $a>b$ be objects of indices $p$ and $q$ respectively, and $\cc^a_b \subset \cc$ the full subcategory endowed with a framed embedding  as in Definition \ref{framedembed}.  This involves a choice of integer $L$ such that for each $\alpha > \beta$ objects in $\cc^a_b$, there is a framed embedding of manifolds with corners, 
 $$
 e_{\alpha, \beta} : Mor (\alpha, \beta) \subset  Mor (\alpha, \beta) \times \br^L\hk \br^L \times J(\mu (\alpha), \mu (\beta)).
 $$
 We now define the functor, $Z^{a,b}_{\cc} : \calj^p_q \to  Sp_*$.     For an integer $m$ with $q \leq m \leq p$,
define 
\begin{equation}\label{zm}
Z^{a,b}_{\cc}(m) = \Sigma^{L-m}\bigvee_{\mu(\alpha) = m}S^{m} _\alpha
\end{equation} where the wedge is taken over all objects $\alpha$ with $a \geq \alpha \geq b$ and $\mu(\alpha) = m$.
Here   the spheres in this notation actually are representing the suspension spectra of the spheres indicated. We continue this abuse of notation whenever dealing with functors with values in $Sp_*$.  The Pontrjagin-Thom construction on the embedding $e_{\alpha, \beta}$
defines a map of one point compactifications, 
\begin{align}\label{pthom}
Z_{\cc}(\alpha, \beta) :  \left(\br^L \times J(\mu (\alpha), \mu (\beta))\right)\cup \infty  &\la  \left(Mor (\alpha, \beta) \times \br^L\right) \cup \infty  \xr{proj} S^L \notag \\
S^L \wedge J(\mu(\alpha), \mu (\beta))^+ &\la S^L \notag 
\end{align}
which we think of as a map of suspension spectra
\begin{equation}\label{zalbe}
Z_{\cc}(\alpha, \beta):  \Sigma^{L-\mu(\alpha)} S^{\mu(\alpha)}_\alpha \wedge J(\mu(\alpha), \mu (\beta))^+ \la \Sigma^{L-\mu(\beta)} S^{\mu(\beta)}_\beta.
\end{equation}
Wedges of these maps (in the category $Sp_*$)  define the functor $Z^{a,b}_{\cc}$ on morphisms. 
 The compatibility conditions of these framings given in Definition \ref{framedembed} assure that the appropriate compositions are preserved. 

\end{proof}

\begin{theorem}\label{celldecomp}
 Let $\cc$ be a compact, smooth, framed  Morse-Smale category  of finite type. Let  $a > b$ be objects.  Then the geometric realization of the functor,
 $$
 Z^{a,b}_{\cc} : \calj^p_q \to  Sp_* 
 $$
 has a cell decomposition with one cell for every object $\alpha$ with $a \geq \alpha \geq b$.  Furthermore, its cellular chain complex,
 $$
 \to \cdots C_m^{a,b} \xr{\p} C_{m-1}^{a,b} \to \cdots  
 $$
  has      boundary homomorphisms  that can be computed in the following way. If $\alpha$ and $\beta$ are objects with $\mu (\alpha) = m $ and $ \mu (\beta) =m-1$,  then the coefficient  
$$
\langle \p(\alpha), \beta \rangle = n_{\alpha, \beta} \in \bz
$$
is given by   $$n_{\alpha, \beta} =   \#\cm(\alpha, \beta).$$
This  number is counted with sign determined  by the framing.  
\end{theorem}

\begin{proof}

Recall    the decomposition (\ref{cell})  into iterated cones of the geometric realization $|Z|$ of  a functor $  Z : \calj^p_q \to  Sp_*$.  In the case of the functor $Z^{a,b}_{\cc} : \calj^p_q \to  Sp_*$ coming from a compact, smooth, Morse-Smale category $\cc$ of finite type, each $Z^{a,b}_{\cc}(m)$ is a wedge of spheres  indexed on the objects, which   are the critical points.  Thus this decomposition gives a cell decomposition of $Z^{a,b}_{\cc}$ with one cell for every critical point. The 
dimension of a cell corresponding to  a critical point   $\alpha$, is $\mu(\alpha) + (L-q).$      Furthermore, the corresponding chain complex (\ref{complex})  is then generated by the objects with $a > \alpha > b$.  The boundary homomorphism in this chain complex is determined by the attaching maps in the cell decomposition (\ref{cell}). It can be computed by observing that if $\alpha$ and $\beta$ are objects with $ \mu (\alpha) =m$ and $\mu (\beta) = m-1$,  then the coefficient  of the boundary homomorphism,   
$$
\langle \p(\alpha), \beta \rangle = n_{\alpha, \beta} \in \bz
$$
 is the degree of the attaching map,
\begin{align}
\phi_{\alpha, \beta} : \Sigma^{L-q-1}S^{m+1}_{\alpha} \xr{\iota_\alpha} |Z^{a,b}_{\cc}|^{(m+L-q )} &\xr{project}  |Z^{a,b}_{\cc}|^{(m+L-q)}/ |Z^{a,b}_{\cc}|^{(m+L-q-1)} \notag \\
&= \Sigma^{L-q}\bigvee_{\mu(\gamma) = m}S^{m} _\gamma \xr{p_{\beta}} \Sigma^{L-q}S^m_\beta. \notag
\end{align}
Here $|Z^{a,b}_{\cc}|^{(r)}$ is the $r^{th}$ skeleton of $|Z^{a,b}_{\cc}|$,    $\iota_{\alpha}$ is the  attaching map of the cell corresponding
to $\alpha$, and $p_{\beta}$ is the projection onto the sphere corresponding to $\beta$.  
  But by definition, this map is constructed as the Pontrjagin-Thom construction on the framed, compact zero dimensional manifold $\cm(\alpha, \beta)$.  In other words, the degree $$n_{\alpha, \beta} = \#\phi_{\alpha, \beta}^{-1}(\infty) = \#\cm(\alpha, \beta)$$
  where this number is counted with sign determined  by the framing.   
\end{proof}

  As described in \cite{cjs}, the standard example of a compact framed Morse-Smale category is the \sl flow category \rm of a Morse function on a closed manifold,
  $$
  f : M \to \br
  $$
  satisfying the Morse-Smale transversality condition. This category, which we denote
  by $\cc_f$  has the set of critical points of $f$ as its objects.  The space of morphisms between critical points $a$ and $b$  is the compact space of piecewise gradient flow lines connecting $a$ to $b$,  $\bcm (a,b)$.   This category was first defined in \cite{cjscat} and studied in \cite{cjs}.     It is well known that the space  of flows $\bcm (a,b)$ is compact framed manifold 
  with corners of dimension $ind (a) - ind (b) -1$.   

\section{The Floer theory}

Our goal in this section is to show that the flow category of the (perturbed) symplectic action functional on the loop space of the cotangent bundle, $L(T^*M)$, is a compact, smooth, framed  Morse-Smale category  of finite type,  assuming one chooses the Hamiltonian and the almost complex structure appropriately.  By the results of the last section
this defines a functor $Z : \calj \to Sp_*$, and therefore a ``Floer homotopy type" which has one cell for every critical point (periodic orbit of the Hamiltonian). This will prove Theorems \ref{frameable} and \ref{floerhomotopy} as stated in the introduction.  

\med

We begin by recalling the basic analytic setup for the Floer theory of $T^*M$, as   described by Abbondandolo and Schwarz \cite{abboschwarz}.  

\subsection{Analytic setup}

Let $M$ be a connected, closed, orientable manifold of dimension $n$, and $T^*M$ its cotangent bundle with its canonical symplectic form $\omega$.   Recall that $\omega$ is exact,   $\omega = d \theta$, where $\theta$ is the  Liouville $1$-form on $T^*M$, where for $x \in M$, and $t\in T^*_xM$, $\theta (x,t)$ is the given by the composition,
 $$
 \theta (x,t): T_{x,t}(T^*M) \xr{dp}T_xM \xr{t}\br.
 $$
 
  Let 
 $$
 H : \br/\bz \times T^*M \to \br
 $$
 be a $1$-periodic Hamiltonian, with corresponding Hamiltonian vector field $X_H$ defined by
 $$
 \omega (X_H(t,x), v) = -dH_{(t,x)} (v)
 $$
 for all $(t,x) \in T^*M$, and $v \in T_{t,x}(T^*M)$.  Here $x \in M$, and $t \in T_x^*M$.
 We will be considering the space of $1$-periodic solutions, $\cp (H)$, of the Hamiltonian equation
 $$
 \frac{dx}{dt} =  X_H(t, x(t)).
 $$
 
 As in \cite{abboschwarz} we make the nondegeneracy assumption,
 
 \bfl
 \bf (H0) \rm Every solution $a \in \cp (H) $ is nondegenerate.  This means that if $\phi^t_H$ is the integral flow of the vector field $X_H$, then $1$ is not an eigenvalue of $D\phi^1_H(x(0) \in GL(T_{x(0)}T^*M)$.  
 
 \efl
 This condition is known to hold for a generic set of $H$ (see \cite{Salzehn}, \cite{weber}).
 
 \med
 Continuing to follow \cite{abboschwarz},  let $L(T^*M)$ denote the space of all loops $x : S^1 \to T^*M$ of Sobolev class $W^{1,2}$.  $L(T^*M)$ has a canonical Hilbert manifold structure.  Let $H$ be a Hamiltonian satisfying condition (H0).  Consider the perturbed symplectic action functional 
 \begin{align}
 \ca_H : L (T^*M)   &\to \br \notag \\
 x &\to \int x^*(\theta -Hdt) = \int_0^1(\theta (\frac{dx}{dt}) - H(t, x(t)) dt.
\end{align}
This is a smooth functional, and its   critical points  are the elements of $\cp (H)$.  
Now let $J$ be a $1$-periodic, smooth almost complex structure on $T^*M$, so that for each $t \in \br/\bz$, 
$$
\langle \zeta, \xi \rangle_{J_t} = \omega (\zeta, J(t,x)\xi ), \quad \zeta, \xi \in T_{x}T^*M, \, x \in T^*M,
$$
is a loop of Riemannian metrics on $T^*M$.  One can then consider the gradient of $\ca_H$ with respect to the metric, $\langle \cdot , \cdot \rangle$, written as
$$
\nabla_J\ca_H(x) = -J(x,t)(\frac{dx}{dt}-X_H(t,x)).
$$
The (negative) gradient flow equation on a smooth curve $u : \br \to L(T^*M)$,
$$
\frac{du}{ds} + \nabla_J\ca_H(u(s))
$$
can be rewritten as a perturbed Cauchy-Riemann PDE, if we view $u$ as a smooth map
$\br/\bz \times \br \to T^*M$,  with coordinates, $t \in  \br/\bz, \, s \in \br$,
\begin{equation}\label{cauchy}
\p_s u -J(t,u(t,s))(\p_t u - X_H(t, u(t,s)) = 0.
\end{equation}

Let $a, b \in \cp (H)$. Abbondandolo and Schwarz defined the space of solutions
 
\begin{align}\label{wab}
W(a,b; H,J) = \{u : \br &\to L(T^*M)\, \text{ a solution to} \,  (\ref{cauchy}), \, \text{such that}  \\
\lim_{s\to -\infty}u(s) &= a, \, \text{and} \,  \lim_{s\to +\infty}u(s) = b \}.\notag
\end{align}
 We call this space $W(a,b)$ for short.  Note that our notation differs from that of \cite{abboschwarz}.  They use the notation $\cm (a,b)$.  We will reserve this notation for the ``moduli space" obtained by dividing out by the free $\br$-action,
\begin{equation}\label{flows}
\cm (a,b) = W(a,b)/\br.
\end{equation}

Give $M$ a Riemannian metric. Let $\nabla$ represent the corresponding Levi-Civita connection.   Abbondandolo and Schwarz then imposed the following further conditions on the Hamiltonian:
\bfl
\bf (H1)\rm There exist $h_0 > 0$  and $h_1 \geq 0$ such that
$$
dH(t, q, p)[\eta] -H(t, q, p) \geq h_0|p|^2 - h_1,
$$
for every $(t, q, p) \in \br/\bz \times T^*M$ (so that $q \in M$, and $p \in T^*_qM$.).

\med
\bf (H2) \rm There exists $h_2 \geq 0$ such that 
$$
|\nabla_q H(t, q, p)|\leq h_2(1 + |p|^2),   \quad |\nabla_q\p H(t, q, p)| \leq h_2(1 + |p|),
$$
for every $(t, q, p) \in \br/\bz \times T^*M$. 

\efl
 As observed in \cite{abboschwarz},  Condition (H1) is a condition of quadratic growth and infinity, and neither conditions (H1) nor (H2) depend on the choice of metric.
 
 An important property of Hamiltonians that satisfy conditions (H0) - (H2) is the following.
 
 \med
 \begin{lemma}(\cite{abboschwarz}) \label{finite}.  Suppose $H$ is a Hamiltionian on $T^*M$ satisfying conditions (H0), (H1), and (H2).  Then for any real number $r$, the set of solutions $a \in \cp (H)$ with $\ca_H(a) \leq r$ is finite.
 \end{lemma}
 
 \med
The Levi-Civita connection coming from  the Riemannian structure on $M$ defines a splitting of the tangent bundle of $T^*M$,
\begin{equation}\label{split}
T(T^*M) \cong p^*(TM) \oplus p^*(T^*M),  
\end{equation}
where $p : T^*M \to M$ is the projection.  This determines a canonical almost complex structure $\hat J$ on $T^*M$ compatible with the symplectic structure, which, with respect to this splitting is given by
$$
\hat J = \begin{pmatrix} 0 & I \\
-I & 0 \end{pmatrix}
$$
In \cite{abboschwarz} the   the spaces $W(a,b)$ were described as the zeros of an appropriate vector field defined by the perturbed Cauchy-Riemann operator (\ref{cauchy}), as follows. 

Define the space $\cb(a,b)$ to be the space of maps $u : \br \times \br/\bz \to T^*M$ of Sobolev class $W^{1,r}_{loc}$ such that there is an $s_0 \geq 0$ for which
$$
u(s,t) = \begin{cases}  exp_{a(t)}(\zeta^-(s,t)) \quad \text{for} \quad s \leq -s_0 \\
exp_{b(t)}(\zeta^+(s,t)) \quad \text{for} \quad  s \geq s_0,  \end{cases}
$$
where $\zeta^-$ and $\zeta^+$ are $W^{1,r}$ sections of $a^*(TT^*M)$ and $b^*(TT^*M)$ respectively.  These pullback bundles are viewed as bundles over $(-\infty , -s_0] \times \br/\bz$ and $[s_0, +\infty) \times \br/\bz$ respectively.        $\cb(a,b)$ has the structure of a Banach manifold, and the tangent space at $u \in \cb(a,b)$ is the space of $W^{1,r}$-sections of $u^*(TT^*M)$.     Notice that there is a homotopy equivalence,
\begin{equation}\label{paths}
\cb(a,b) \simeq \Omega_{a,b}(L(T^*M))  
\end{equation}
where $\Omega_{a,b}(L(T^*M))$ is the space of continuous paths $\gamma : [0,1] \to L(T^*M)$ with
$\gamma (0) = a$, and $\gamma (1) = b$.  This space has the compact-open topology.

Define $\ct \cb \to \cb(a,b)$ to be the Banach - bundle whose fiber $\ct_u (\cb (a,b))$ at $u \in \cb (a,b)$ is the space of $L^r$-section of $u^*(TT^*M)$.  Then $W(a,b)$ is the space of zeros of the smooth section,
\begin{align}\label{vfield}
\p_{J,H} :  \cb (a,b) &\to \ct \cb(a,b)   \\
u &\to \p_s u + \nabla_J\ca_H(u) = \p_s u - J(t,u)(\p_tu - X_H(t,u(s,t)). \notag
\end{align}

Define $\cj_{reg}(H)$ to be the set of all time dependent, periodic smooth $\omega$-compatible almost complex structures such that $\|J - \hat J\| < \infty$   and such that the section $\p_{J,H}$ is transverse to the zero section, for every $a, b \in \cp_H$.     They observe that this is a residual subspace of the space of all $\omega$-compatible almost complex structures $J$ with   $\|J - \hat J\| < \infty$.   Thus the following condition is generic.

\med
 \bfl
 \bf (J1). \rm   $J \in \cj_{reg}(H)$.
 
 \efl

\med
\begin{theorem}(\cite{abboschwarz})\label{smooth}  Assume the Hamiltonian $H$   satisfies conditions
(H0), (H1), (H2).      Then there exists a number $j_0 > 0$ such that if  a time dependent  $\omega$-compatible almost complex structure $J$ satisfies condition (J1) with respect to the Hamiltonian $H$,  and   $\|J - \hat J\| < j_0$, then the spaces $W(a,b; H, J)$ are  all precompact,
 orientable smooth manifolds of dimension $\mu (a) - \mu (b)$, where $\mu (x)$ is the Conley-Zehnder index of the periodic solution $x \in \cp (H)$.  
\end{theorem}

\med
\noindent
\bf Remark.  \rm The compactness results in the above theorem basically follow from the fact that the 
  symplectic form $\omega$ on $T^*M$ is exact, and hence there can be no bubbling.  See \cite{abboschwarz} for details.

\med

If $\cm (a,b) = W(a, b)/\br$ is the  moduli space, then construct the space  $\bar \cm(a,b)$   of  ``piecewise flow lines" in the usual way:
$$
\bar \cm (a,b) = \bigcup_{a=a_1 >a_2> \cdots > a_k = b} \cm(a_1, a_2) \times \cdots \times \cm (a_{k-1}, a_k),
$$
Here the partial order is given by $a_1 \geq a_2$ if $W(a_1, a_2) \neq \emptyset$.

The topology of $\bar \cm (a,b)$ can be described as follows.  
Since the action functional $\ca_H$ is strictly decreasing along flow lines,  an element $u \in  W(a,b)$  determines a diffeomorphism $\br \cong (\ca_H(b), \ca_H(a))$ given by the composition,
$$
 \br  \xr{u} L(T^*M) \xr{\ca_H}  \br.
 $$
 This defines a parameterization of any $\gamma \in \cm(a,b)$ as a map
 $$
\gamma : [ \ca_H(b), \ca_H(a)] \to L(T^*M)
$$
that satisfies the differential equation 
\begin{equation}\label{diffeq}
\frac{d\gamma}{ds} =  \frac{\nabla \ca_H (\gamma (s))}{|\ca_H(\gamma (s))|^2},
\end{equation}
as well as the boundary conditions \begin{equation}\label{bound}
\gamma (\ca_H(b)) = b  \quad \text{and} \quad \gamma (\ca_H(a))=a.
\end{equation}  From this viewpoint, the \sl compactification \rm $\bcm(a,b)$ can be described as the space of all continuous maps $ [ \ca_H(b), \ca_H(a)] \to L(T^*M)$ that are piecewise
smooth, (and indeed smooth off of the critical values of $\ca_H$ that lie between $\ca_H(b)$ and $\ca_H(a)$), that satisfy the differential equation (\ref{diffeq}) subject to the boundary conditions (\ref{bound}).  It is topologized with the compact open topology.

We now invoke Gluing Assumption (\ref{glue}), which   implies  that the  spaces $\bar \cm(a,b)$ are compact, smooth, orientable manfolds with corners.
The dimension of $\bar \cm (a,b)$ is $\mu (a) - \mu(b) -1$.   For the rest of the paper we operate under this gluing assumption.
The following is immediate.

\med
\begin{proposition}.  The compact moduli spaces $\bar \cm (a,b)$ are $\langle  k(a,b)\rangle$ - manifolds, where $k(a,b) = \mu (a) - \mu (b) -1$.  
\end{proposition}

\med
\begin{proof}  The faces of $\bar \cm (a,b)$ are given as follows.  For $j = 1, \cdots , k(a,b)$,  let 
$$
F_j (\bar \cm (a,b)) = \bigcup_{\mu (c) = \mu(a) - j} \bar \cm (a, c) \times \bar \cm (c,b).
$$
This face structure clearly satisfies the intersection property necessary for being a $\langle k(a,b) \rangle$- manifold (See  the definition   in the beginning of section 1.2.)
\end{proof}

\med
The following is also straightforward.

\med
\begin{proposition}\label{htpyequiv} The inclusion of the moduli space into its compactification,
 $$
 \iota: \cm(a,b) \hk \bcm (a,b)
 $$
 is a homotopy equivalence.
 \end{proposition}
 \begin{proof}  By the gluing assumption (\ref{glue}), one has local diffeomorphisms,
 $\cm(a,a_1) \times \cdots \times \cm(a_q,b) \times [0, \eps) \to \bcm (a,b)$, for every sequence of intermediate critical points, $a \geq a_1 \geq a_2 \geq \cdots \geq a_q \geq b$.  The boundary $\p \bcm(a,b)$ consists of the images of these gluing maps restricted to
 $\cm(a,a_1) \times \cdots \times \cm(a_q,b) \times \{0\}$, and these gluing maps restrict on the open intervals to give local diffeomorphisms,  $\cm(a,a_1) \times \cdots \times \cm(a_q,b) \times (0, \eps) \to \cm (a,b).$ 
 Thus 
 $$
 \cm(a,b) = \bcm(a,b) \, - \,  \bigcup_{a \geq a_1 \geq a_2 \geq \cdots \geq a_q \geq b} 
 \cm(a,a_1) \times \cdots \times \cm(a_q,b) \times \{0\}
 $$
 where the union is taken over all such sequences of intermediate critical points.  In this notation we are identifying $\cm(a,a_1) \times \cdots \times \cm(a_q,b) \times \{0\}$ with its image in $\bcm(a,b)$. 
 
 Let $\ck (a,b)$ be the complement of a smaller open collar around the boundary:
 $$
 \ck(a,b) = \bcm(a,b) \, - \,  \bigcup_{a \geq a_1 \geq a_2 \geq \cdots \geq a_q \geq b} 
 \cm(a,a_1) \times \cdots \times \cm(a_q,b) \times [0, \eps/2)
 $$
 In particular, $\ck(a,b)$ is a subspace of the open moduli space $\ck(a,b) \subset \cm(a,b).$
 But clearly $\ck(a,b)$ is homeomorphic to $\bcm(a,b)$.  Such a homeomorphism is induced by the natural affine homeomorphism between the intervals $[\eps/2, \eps) $ and $[0, \eps)$.
 With respect to this homeomorphism the inclusion $\ck(a,b) \subset \cm(a,b)$ is a clearly a homotopy inverse to the inclusion $\cm(a,b) \subset \bcm(a,b)$.    
  \end{proof}
  
  The following says that the stable normal bundle of $\bcm(a,b)$ is, up to isomorphism, completely determined by the stable normal bundle of the open moduli space $\cm (a,b)$.  In particular it does not depend on the actual gluing maps in gluing assumption (\ref{glue}). 
  
  \med
  \begin{corollary}\label{stabnormal} The stable normal bundle of $\bcm(a,b)$ is classified by any map $\bcm(a,b) \to BO$ whose restriction to the open moduli
  space $\cm (a,b)$ classifies its stable normal bundle.
  \end{corollary}

 \subsection{A framing of the moduli space of $J$-holomorphic cylinders}
  
  The goal of this section is to prove Theorem \ref{frameable} in the introduction, and to define the Floer homotopy type of $T^*M$.   We now restate this Theorem \ref{frameable}  with the appropriate hypotheses.
  
  \med
  
  \begin{theorem}\label{trivial}  Let  $H$ and $J$  satisfy the conditions of Theorem \ref{smooth}. Then the compactified moduli spaces, 
$\bcm (a,b) = \bcm(a,b; H, J) $ have  framings on their stable normal bundles. Here framing means as manifolds with corners as in Definition (\ref{framedcorner}).  
\end{theorem}

  \med

  \begin{proof}  In \cite{cjs} the authors discussed how the basic obstruction to the framability  of the morphism spaces $Mor (a,b)$ in the flow category $\cc_N$ of the
symplectic action functional on the loop space, $LN$, of a symplectic manifold, $(N^{2n}, \omega)$, is the ``polarization class",   which is a homotopy class of map 
$$
\rho : LN \to U/O
$$
where $U/O = \lim_{k\to \infty} U(k)/O(k)$.      The map $\rho$ is defined as follows.  A compatible almost complex structure on $N$ defines a map $\tau: N \to BU(n)$ classifying the isomorphism class of the  tangent bundle as a complex vector bundle.  The homotopy class of $\tau$ does not depend on the particular compatible almost complex structure chosen.  Applying loop spaces, one has a composite map,
\begin{equation}\label{polarization}
\rho: LN \xr{L(\tau)} L(BU(n)) \hk LBU \simeq BU \times U   \to U \to U/O.
\end{equation}
Here the homotopy equivalence $L(BU) \simeq BU \times U$ is well defined up to homotopy, and is given by a trivialization of the   fibration
$$
U \simeq \Omega BU \xr{\iota} L(BU) \xr{ev}BU
$$
where $ev :  LX \to X$ evaluates a loop a $0 \in \br/\bz$.  

\med
The reason we refer to this invariant as the ``polarization class" of the loop space $LN$, is  because when viewed as an infinite dimensional manifold, the tangent bundle $T(LN)$ is polarized, and is classified up to homotopy by the map $\rho$.  See \cite{cjs} for details.

 \med
We now observe that in the case of the cotangent bundle, this polarization class is trivial.

\begin{lemma}\label{canonical}  The polarization class $\rho : L(T^*M) \to U/O$ is has a canonical null homotopy.
\end{lemma}

\begin{proof} As described above (\ref{split}), a Riemannian metric on $M$ defines  via the Levi-Civita connection, a splitting $T(T^*M) \cong p^*(TM) \oplus p^*(T^*M)$, and an $\omega$-compatible almost complex structure defines an isomorphism of the complex
vector bundle $T(T^*M)$, with the complexification of the tangent bundle of $M$,
$$
T(T^*M) \cong p^*(TM)\otimes \bc.
$$
On the level of classifying maps this means the following diagram homotopy commutes:
\begin{equation}
\begin{CD}
M @>\tau_M >>  BO(n) \\
@V\cap Vz V   @VVcV \\
T^*M  @>>\tau_{T^*M}> BU(n)
\end{CD}
\end{equation}
where the horizontal maps $\tau$ classify the relevant tangent bundles, $z : M \hk T^*M$ is the zero section, and $c : BO(n) \to BU(n)$ is the complexification map.  That means we have a homotopy commutative diagram,
$$
\begin{CD}
LM    @>L\tau_M >> LBO(n) @>>> L(BO) \simeq O \times BO  @>project >> O \\
@V\cap Vz V  @VVc V   @VVcV  @VVcV \\
L(T^*M) @>>L\tau_{T^*M}>  LBU(n) @>>>L(BU)\simeq U \times BU @>>project > U @>>> U/O
\end{CD}
$$
Now  the bottom horizontal composition is the polarization class
$\rho : L(T^*M) \to U/O$.  But since the composition $O \xr{c} U \to U/O$ is canonically null homotopic, the commutativity of this diagram says that the composition
$$
LM \xr{z} L(T^*M) \xr{\rho} U/O
$$
is canonically null homotopic.    But the zero section $z : LM \hk L(T^*M) $ is a homotopy equivalence, with homotopy inverse $Lp: L(T^*M) \to LM$.  This implies that
$\rho: L(T^*M) \to U/O$ is canonically null homotopic.
\end{proof}

 In order to complete the proof of Theorem \ref{trivial},  we will now   show how  the  canonical null homotopy of the polarization class  induces canonical trivializations of the stable normal bundles of the morphism manifolds, $\bcm (a,b)$. This relationship was alluded to in \cite{cjs}.

 Consider the space $W(a,b)= W(a,b; J,H)$ described above.  We first show that under the hypotheses of this theorem, $W(a,b)$ has a stably trivial tangent bundle.  Recall that the manifold $W(a,b)$ can be described as the space of zero's of a vector field  (\ref{vfield}).  The  tangent  space of $u \in W(a,b)$ is therefore given by the kernel of the fiberwise derivative,  
$$
D_{fib}\p_{J,H}(u) : T_u \cb(a,b) \to \ct_u\cb (a,b).
$$   A trivialization $\Phi$ of $u^*(TT^*M)$ defines a conjugacy between $D_{fib}\p_{J,H}(u)$ and a bounded operator, 
$ 
D_S : W^{1,r}(\br \times \br/\bz; \br^{2n}) \la L^r(\br \times \br/\bz; \br^{2n})$
of the form
\begin{equation}\label{local}
D_S(v) = \p_s(v) - J_0\p_t(v) - S(s,t)v
\end{equation}
where $S$ is a smooth family of endomorphisms of $\br^{2n}$ in which the limits
$$
S^{\pm}(t) = \lim_{s\to \pm \infty} S(s,t)
$$
are symmetric.

Thus the stable tangent bundle of $W(a,b)$  is classified in the following way.
Consider the fiber bundle
$$
Fr  \xr{p} W(a,b)
$$
whose fiber over $u \in W(a,b)$ is the space of   Fredholm operators, $T_u \cb(a,b) \to \ct_u\cb (a,b)$.    Notice that this bundle is trivial, because it is induced from   principal $GL(W^{1,r}(\br \times \br/\bz; \br^{2n})$ and $GL( L^r(\br \times \br/\bz; \br^{2n})$-bundles.   But by Kuiper's theorem, these general linear groups   are contractible.  This says that the bundle, $
Fr  \xr{p} W(a,b)$ is trivial, and  the space of trivializations is contractible.  
Thus the assignment to $u \in W(a,b)$, the operator $D_{fib}\p_{J,H}(u)$ is a section
of $Fr  \xr{p} W(a,b)$, which gives a well defined homotopy class of map to the space of Fredholm operators,
$$
D_{fib}\p_{J,H} : W(a,b) \la Fred (W^{1,r}(\br \times \br/\bz; \br^{2n}) ,  L^r(\br \times \br/\bz; \br^{2n}).
$$
This space of Fredholm operators is homotopy equivalent, via Atiyah's theorem \cite{atiyah} to the classifying space $\bz \times BO$.  Since these operators all have index
equal to the difference of the Conley-Zehnder indices, $\mu (a) - \mu (b)$  the image of this map lies in the component of $\bz \times BO$ corresponding to $\mu (a) - \mu (b) \in \bz$.   By the above description of the tangent spaces $T_u(W(a,b))$, the induced (homtopy class of) map
\begin{equation}\label{stabletan}
D_{fib}\p_{J,H} : W(a,b) \to BO
\end{equation}
classifies the stable tangent bundle.  We therefore must show that this map is null homotopic. 

\med
Notice that since $D_{fib}\p_{J,H}(u)$ is a perturbation of the Cauchy-Riemann operator,
the following diagram homotopy commutes:

\begin{equation}\label{diagram}
\begin{CD}
W(a,b)   &   @>D_{fib}\p_{J,H} >>    &  \bz \times BO \\
@V\cap VV &&  @AA \iota A \\
\cb (a,b)  \simeq \Omega_{a,b}(L(T^*M))  @>>\Omega L (\tau )>  \Omega L BU(n)  @>>ind \, \bar \p >  \bz \times BU  
\end{CD}
\end{equation}
   The map $ind \, \bar \p$ assigns to a map
$u : S^1 \times \br \to BU(n)$ the index of the $\bar \p$-operator coupled to $u$, as 
constructed by Atiyah in \cite{atiyah2} in his well known proof of Bott-periodicity using elliptic operators. In particular, $ind \, \bar \p$ factors up to homotopy  through a map $ind \, \bar \p  : \Omega^2 BU \xr{\simeq} \bz \times BU$ which is the homotopy equivalence inverse to the Bott map, $\beta : \bz \times BU \to \Omega^2 BU$.  The map $\iota : BU \to BO$ is induced by the inclusion maps, $U(k) \hk O(2k)$.  

In other words,  $ind \, \bar \p :\Omega L BU(n) \to \bz \times BU$ factors as the composition
$$
\begin{CD}
\Omega L BU(n) @>\hk >> \Omega L BU \simeq \Omega^2 BU \times \Omega BU @>proj >> \Omega^2 BU @>ind \, \bar \p >\simeq >  \bz \times BU
\end{CD}
$$
But also notice that standard homotopy theory implies the following diagram homotopy commutes:
$$
\begin{CD}
\Omega^2 BU &  @>ind \, \bar \p >\simeq >  & \bz \times BU \\
@V\simeq VV   &&  @VV\iota V \\
\Omega U @>>project > \Omega U/O @>>\simeq > \bz  \times BO.
\end{CD}
$$

Putting these facts together with the definition of the polarization class  
$\rho : L(T^*M) \to U/O$,  we have the following homotopy commutative diagram:
\begin{equation}\label{tangent}
\begin{CD}
W(a,b)      @>D_{fib}\p_{J,H} >>     \bz \times BO \\
@VVV   @VV \simeq V \\
\Omega L(T^*M) @>>\Omega \rho >  \Omega (U/O)
\end{CD}
\end{equation}

Thus the canonical  null homotopy of the polarization class, $\rho: L(T^*M) \to U/O$ (Lemma \ref{canonical})  defines
a null homotopy of the induced map of loop spaces, $\Omega \rho : \Omega L(T^*M) \to \Omega (U/O)$, which, by the above diagram defines null homotopies of the classification
maps of the stable tangent bundles, $D_{fib}\p_{J,H}  : W(a,b) \to \bz \times BO$.

Now the spaces of flows, $\cm (a,b)$ can be viewed as a subspace of $W(a,b)$ as the intersection of $W(a,b)$ with a level set of the action functional $\ca_H : L(T^*M) \to \br$.  
With respect to this embedding, there is a natural isomorphism of bundles,
$$T(\cm(a,b)) \times \br \cong T(W(a,b)).$$
So the stable trivialization of the tangent bundles $TW(a,b)$ induce stable trivializations
fo $T(\cm (a,b))$.  

We need to organize these stable trivializations further in order to  complete the proof of  Theorem \ref{framedcorner}.

Let $a \geq b$ be a pair of critical points of $\ca$.  We define the path space
$\Omega_{a,b}U/O$ to be the space of paths
$$
\theta : [\ca(b), \ca (a)] \to U/O
$$
with the boundary conditions $\theta (\ca (b)) = \rho (b)$,   $ \theta (\ca(a)) = \rho (a)$.
Notice that there is a homotopy equivalence $\Omega_{a,b}U/O \simeq \Omega U/O \simeq \bz \times BO$.   We therefore take $\Omega_{a,b}U/O$ to be our model of $\bz \times BO$, and view it as a constant $\lkr$-diagram, where $k = \mu(a)-\mu(b)-1$.  

Consider the composition
$$
\nu_{a,b}: \bcm(a,b) \hk \Omega_{a,b}LT^*M \xr{\Omega_{a,b}\rho} \Omega_{a,b}U/O \xr{\Omega_{a,b}(-1)}\Omega_{a,b}U/O
$$
where $``-1" : U/O \to U/O$ is a fixed self homotopy equivalence that induces a map of degree $-1$ on $\pi_1(U/O) = \bz$.    By diagram (\ref{tangent}),  the restriction of this map to the open moduli space $\cm(a,b)$ classifies its stable normal bundle, hence
by Corollary \ref{stabnormal}, $\nu_{a,b}$ classifies the stable normal bundle of $\bcm(a,b)$.
Moreover, notice that the following diagrams obviously commute:
\begin{equation}\label{commdiag}
\begin{CD}
\bcm(a,c) \times \bcm(c,b) @>compose >> \bcm(a,b) \\
@V\cap VV  @VV\cap V \\
 \Omega_{a,c}LT^*M \times \Omega_{c,b}LT^*M  @>compose >> \Omega_{a,b}LT^*M \\
  @V  \Omega_{a,c}\rho \times \Omega_{c,b}\rho VV    @VV\Omega_{a,b}\rho V \\
   \Omega_{a,c}U/O \times \Omega_{c,b}U/O     @>compose >> \Omega_{a,b}U/O \\
   @V\Omega_{a,c}(-1)\times \Omega_{c,b}(-1) VV  @VV\Omega_{a,b}(-1) V\\
   \Omega_{a,c}U/O \times \Omega_{c,b}U/O  @>>compose > \Omega_{a,b}U/O  \\
\end{CD}
\end{equation}
This means that the restriction of the map of $\lkr$-diagrams $\nu_{a,b} : \bcm(a,b) \to \Omega_{a,b}U/O$ to a stratum of the boundary (i.e to any object in the poset $\ut^k$)
classifies the stable normal bundle.  Therefore as a map of $\lkr$-diagrams,
$\nu_{a,b}$ classifies the stable normal bundle of the $\lkr$-manifold $\bcm(a,b)$. 

By Definition (\ref{trivial}), a framing of this stable normal bundle map is an appropriate
lifting of $\nu_{a,b}$ to a contractible fibration.  We define this as follows.

Recall that $\pi_1(U/O) \cong \bz$.  Let $f : S^1 \to U/O$ be a generator. That is, $f$ induces an isomorphism on the fundamental group.  Define $\phi: \cs \to U/O$  to be a fibration replacing $f$ in the following way. Let $\cs = \{(t, \gamma) \in  S^1 \times U/O^I \, : f(t) = \gamma (0) \}$.  Here $X^I$ denotes the space of continuous paths, $[0,1] \to X$.  $\phi : \cs \to U/O$ is defined by
$\phi (t, \gamma) = \gamma (1)$.  Notice that $\phi$ is a ``fibrant replacement" of $f$, because $\cs$  has a natural  homotopy equivalence to $S^1$, with respect to which
$\phi$ is homotopic to $f$. 

Notice that the canonical null homotopy of the polarization class $\rho : L(T^*M) \to U/O$ defines a null homotopic lifting,
$$
\trho : L(T^*M) \to \cs
$$
which we now fix.

Like above, let
$\Omega_{a,b}\cs $ to be the space of paths, 
$$
\theta : [\ca(b), \ca (a)] \to \cs
$$
with the boundary conditions $\theta (\ca (b)) = \trho (b)$,   $ \theta (\ca(a)) = \trho (a)$.
Notice that $\Omega_{a,b}\cs  \simeq \Omega \cs \simeq \Omega S^1$, so that it has one path component for every integer, but each path component is contractible.
We notice that $\phi$ induces a fibration, which by abuse of notation we still call
$$
\phi :\Omega_{a,b}\cs \to \Omega_{a,b}U/O
$$
where each path component of $\Omega_{a,b}U/O$ is covered by a contractible
path component of $\Omega_{a,b}\cs.$ Again, we think of this   as a fibration of constant $\lkr$-diagrams, where $k = \mu(a) - \mu (b) -1$.  

The map $\trho : L(T^*M) \to \cs$  then allows us to define the lifting of $\nu_{a,.b}$.
\begin{equation}\label{lift}
\tilde \nu_{a,b} : \bcm(a,b) \hk \Omega_{a,b}L(T^*M) \xr{\Omega_{a,b}\trho} \Omega_{a,b}\cs \xr{\Omega_{a,b}(-1)} \Omega_{a,b}\cs
\end{equation}
where $``-1" : \cs \to \cs$ is a lifting of the degree $-1$-map $-1 : U/O \to U/O$.  
This lifting is again viewed as a map of $\lkr$-diagrams.  Since $\bcm(a,b)$ is connected, the image of $\tilde \nu_{a.,b}$ is in a connected component of $\Omega_{a,b}\cs$, which is a contractible fibration over the corresponding component of $\Omega_{a,b}U/O$, which is our model for $BO$.  This gives a framing of the stable normal bundle of $\bcm(a,b)$ according to Definition (\ref{framedcorner}) above.  \end{proof}

\med
We observe that we have actually proved something more than Theorem \ref{trivial}, which says that the compact moduli spaces have framings.   Namely the analogue of
commutative diagram (\ref{commdiag}) for the framing $\trho$ defined by (\ref{lift}) implies the following.

\begin{corollary}\label{multiply}  The framings $\trho_{a,b} : \bcm \to \Omega_{a,b}\cs$
are multiplicative, in the sense that for any three critical points $a_1 \geq a_2 \geq a_3$, the following diagram of framings commute:
$$
\begin{CD}
\bcm(a_1, a_2) \times \bcm(a_2, a_3) @>compose >> \bcm(a_1, a_3) \\
@V\trho_{a_1,a_2} \times \trho_{a_2, a_3} VV  @VV\trho_{a_1, a_3} V \\
\Omega_{a_1, a_2}\cs \times \Omega_{a_1, a_2}\cs   @>>compose > \Omega_{a_1, a_3}\cs.
\end{CD}
$$
\end{corollary}


\med
In \cite{laures} Laures showed how the  well known relationship between  framings of stable normal    bundles of closed manifolds, and isotopy classes of framed embeddings into high codimension Euclidean space, extends to $\langle k \rangle$ -manifolds.  From this  one immediately sees  that Corollary \ref{multiply} implies the following, which in view of Theorems \ref{fltype} and \ref{celldecomp}  assures us the existence of a ``Floer homotopy type" of   $T^*M$.   

\med
\begin{theorem}\label{chframed}  Let $\cc_H$ be the flow category of the symplectic action
 $\ca_H : L(T^*M) \to \br$  where $H$ is a Hamiltonian and $J$ is an almost complex structure satisfying the conditions of Theorem \ref{smooth} above.    Then $\cc_H$ is a compact, smooth, framed Morse-Smale category of finite type.   In particular for each  pair of critical points  $a>b$, the subcategory $(\cc_H)_b^a$
has a framed embedding.  Furthermore these framed embeddings are compatible in the sense that the induced framing class  of the stable tangent bundle  of $\bcm (\alpha, \beta)$  with $a\geq \alpha > \beta \geq b$ is independent of the choice of $a$ and $b$.   
\end{theorem}

\med
Let us consider the Floer homotopy type  induced by the framed category $\cc_H$, under the assumptions of theorem.  Consider   nonnegative integers,   $p > q  $, and the corresponding subcategories, 
$$(\cc_H)^q \hk (\cc_H)^p \subset \cc_H,$$  where $(\cc_H)^m$ denotes the full subcategory of $\cc_H$ generated by critical points (objects)  less than or equal to $\alpha$.  Recall from \cite{abboschwarz} that for every $m$ there are finitely many critical points with index  $\leq m. $  Thus these are finite categories (i.e have finitely many objects).  They have   framed embeddings,
which induce functors, 
$$
Z_{\cc_H}^{q} : \calj^q_0 \to  Sp_* \quad \text{and} \quad  Z_{\cc_H}^{q} : \calj^p_0 \to  Sp_*
$$
as in Theorem \ref{fltype}.   These  framed embeddings involve choices of integers
$L_q$ and $L_p$ with $L_p \geq L_q$ and  framed embeddings of manifolds with corners,
$$
 e^{q,0}_{\alpha, \beta} : Mor_{(\cc_H)^b_c} (\alpha, \beta) = \bcm(\alpha, \beta) \subset  \bcm (\alpha, \beta) \times \br^{L_q} \hk \br^{L_q} \times J(\mu (\alpha), \mu (\beta)).
 $$ and 
  
  $$
 e^{p,0}_{\alpha, \beta} : Mor_{(\cc_H)^a_c} (\alpha, \beta) = \bcm(\alpha, \beta) \subset  \bcm (\alpha, \beta) \times \br^{L_p} \hk \br^{L_p} \times J(\mu (\alpha), \mu (\beta)).
 $$

 The framings
  of the tangent bundle  of $\bcm (\alpha, \beta)$ coming for the embeddings $ e^{q,0}_{\alpha, \beta} $ and $ e^{p,0}_{\alpha, \beta}$ are compatible.  So the Pontrjagin-Thom constructions  
$$
\tau^{q,0}_{\alpha, \beta} :  S^{L_{q}}\wedge J(\mu (\alpha) , \mu (\beta))^+  \to  \bcm(\alpha, \beta)_+ \wedge S^{L_{q}}     \xr{proj}S^{L_{q}} 
$$
and 
$$
\tau^{p,0}_{\alpha, \beta} :  S^{L_{p}}\wedge J(\mu (\alpha) , \mu (\beta))^+  \to  \bcm(\alpha, \beta)_+ \wedge  S^{L_{p}}      \xr{proj} S^{L_{p}} 
$$
have the property that  $\tau^{p,0}_{\alpha, \beta}$ is the $(L_{p} - L_{q})$-fold suspension of     $\tau^{q,0}_{\alpha, \beta}$.  This defines a map on geometric realizations,

$$
\eps_{a,b,c} : \Sigma^{(L_{p} - L_{q})} |Z_{\cc_H}^{q}| \to |Z_{\cc_H}^{p}|.
$$
This structure defines a spectrum (up to homotopy), which we abbreviate by $Z(T^*M)$. 
This is the Floer homotopy type of $T^*M$.   

\med
We note that, as described in \cite{cjs}, in general the Floer homotopy type  is a prospectrum.  What allowed us to obtain a spectrum in this situation was that the critical points were bounded below under the partial ordering.  This allowed us to consider the framings associated to the categories, $(\cc_H)^p$.  When one is dealing more generally with a compact, framed category $\cc$ that does  not have this boundedness condition, one has to take an inverse system of the geometric realizations
of the functors $Z^{a,b}_{\cc}: J^{\mu(a)}_{\mu(b)} \to  Sp_*$ as the index   $\mu(b)$ decreases.  This produces the inverse system of spectra making up the prospectrum.
  In this case, however, the Floer homotopy type of $T^*M$ is given by the spectrum $Z(T^*M)$.  By Theorems \ref{fltype} and \ref{celldecomp} we have the following consequence, which proves Theorem \ref{floerhomotopy}  in the introduction.

\med
\begin{corollary}\label{floerspec} Given a Hamiltonian $H$ and almost complex structure $J$ satisfying the hypotheses of Theorem \ref{smooth}, the Floer homotopy type of  $T^*M$  is the $C.W$- spectrum $Z(T^*M)$ which has one cell for each periodic solution $a \in \cp (H)$.
The corresponding cellular chain complex is the Floer complex, and hence
$$
H_*(Z(T^*M)) \cong HF_*(T^*M).
$$
\end{corollary}

\section{The Floer homotopy type of $T^*M$ and the free loop space.}
By Corollary \ref{floerspec} and the results of Abbondandolo and Schwarz \cite{abboschwarz}, we know that the Floer homotopy type $Z(T^*M)$ has homology isomorphic to that of the free loop space, $H_*(Z(T^*M)) \cong H_*(LM)$.  The goal of this section is to strengthen this by proving Theorem \ref{equivalence} as stated in the introduction.     Our method in comparing the homotopy types  of  $Z(T^*M)$ and $\Sigma^\infty(LM_+)$ is to adapt the techniques of Abbondandolo and Schwarz \cite{abboschwarz}  that produced an isomorphism between  the Morse chain complex from a specific Palais-Smale Morse function on $LM$, with the Floer chain complex of $T^*M$ associated to  $\ca_H$.  We will adapt their techniques
to compare the framed bordism type of the moduli spaces of flow lines in each case, in order to conclude a relationship between the corresponding stable homotopy types.
 
 {\it Proof of Theorem \ref{equivalence}.}  We begin by recalling  the Morse theory on the loop space, $LM$ coming from a Lagrangian dynamical system, as studied in \cite{abboschwarz}.   In this setting there is a smooth Lagrangian
$$
L : \br/\bz \times TM \to \br
$$
satisfying and appropriate strong convexity property, bounds on its second derivatives, as well as   nondegeneracy properties (see properties (L0) (L1), and (L2) as defined in section 2 of \cite{abboschwarz}).    In such a setting,
the Hamiltonian 
\begin{align}
H : \br/\bz \times T^*M &\to \br  \notag \\
H(t, q,p) &=  \max_ {v \in T_qM} (p(v) - L(q,t,v))
\end{align}
satisfies conditions  (H0), H(1), and (H2) above, which are the conditions that we have assumed throughout.  Moreover in this context, the energy functional,
\begin{align}
\ce : LM &\la \br   \notag \\
\ce (\gamma) &= \int_0^1 L(t, \gamma (t), \frac{d \gamma}{ds}(t)) dt  \notag
\end{align}
is a Morse function that is bounded below,  its critical points, $\cp_L$,  have finite Morse indices,  and $\ce$ satisfies the Palais-Smale condition.   Here $M$ is given a Riemannian metric, and   $LM$ is   an appropriate  compatible Hilbert manifold model for the loop space.   The Legendre transform then gives a bijective correspondence between the critical points of $\ce$, $\cp_L$, and the critical points of $\cah$,  $\cp (H)$,
defined by 
$$
\cl (c) = (c, g(\frac{dc}{dt}))
$$
where $g : TM \xr{\cong} T^*M$ is the isomorphism induced by the metric.  Indeed
$\cl$ can be viewed as  a homotopy equivalence,
\begin{align}\label{legendre}
\cl : LM &\xr{\simeq} L(T^*M)) \\
\gamma &\to (\gamma, g (\frac{d\gamma}{dt})) \notag
\end{align}
It was shown in \cite{abboschwarz} that for any loop $(\gamma, \eta)  \in L(T^*M)$, then
\begin{equation}\label{inequal}
\cah (\gamma, \eta) \leq \ce (\gamma)
 \end{equation}
with equality holding only when $\eta (s) \in T^*_{\gamma (s)}M$ is given by
$\eta (s) = g (\frac{d\gamma}{dt}(s))$, that is, $(\gamma, \eta) = \cl (\gamma).$

The  Morse complex for the function $\ce : LM \to \br$, which we call $C_*^{\ce}(LM)$,  is generated by $\cp_L$, where the Floer complex, $CF_*(T^*M)$ is generated by $\cp (H)$.  However the Legendre transform correspondence of these generating sets does not yield a chain map.  Nevertheless, Abbondondolo and Schwarz found and isomorphism of chain complexes that we call
  \begin{equation}\label{asiso}
\begin{CD}
\Psi_* : C_*^{\ce}(LM)  @>\cong >> CF_*(T^*M),
 \end{CD}
 \end{equation}
 which yields an isomorphism, $H_*(LM) \cong HF_*(T^*M)$.    $\Psi_*$ is defined by studying zero and one dimensional moduli spaces of  ``mixed flow lines",  $\cm^+(a, \beta)$, where $a \in \cp_L$ and $\beta \in \cp (H)$.  Adapting their ideas, we will use the higher dimensional moduli spaces, to show that this chain isomorphism is induced by a   equivalence of spectra,
 $$
\Psi : \Sigma^\infty (LM_+) \xr{\cong} Z(T^*M),
 $$
 thus proving Theorem \ref{equivalence}.   We now recall the definition
 of the moduli spaces $\cm^+(a, \beta)$ from \cite{abboschwarz}, as well as some of their basic properties.  In \cite{abboschwarz} it was proved that these moduli spaces have orientations.  The main technical result we need to prove Theorem \ref{equivalence}
 is that in fact they have natural framings, that extend the framings of the spaces $\bcm (\alpha, \beta; J, H)$ studied above.  As before, let $p : T^*M \to M$ be the projection map. 
 
 \begin{definition}(\cite{abboschwarz} section 3)  For $a \in \cp_L$ and $\beta \in \cp (H)$,  let $\cm^+(a, \beta)$ be  the set of all $C^\infty$ maps 
 $$
 u : \br_+ \times S^1 \to T^*M
 $$
 of Sobolev type $W^{1,r}$, such that 
 \begin{align}
\p_su - J(t,u)(\p_tu - X_H(t,u)) &=0 \quad \text{on} \, (0, \infty) \times S^1 \notag \\
p \circ u (0, \cdot) \in W^u(a) &= W^u(a; \ce), \quad \text{the unstable manifold of the critical point $a$}   \notag \\
\lim_{s \to +\infty} u(s,t) &= \beta (t) \quad \text{uniformly in } \, t \in S^1. \notag
\end{align}
\end{definition}

It was proved in \cite{abboschwarz} that the spaces, $\cm^+(a,\beta )$ are, (with the choices of $J$ and $H$ as above),   smooth, oriented manifolds of dimension $ ind (a) - \mu(\beta)$, where $ind (a)$ is the Morse index.  They showed that these manifolds are precompact, and have compactifcations, $\bcm^+(a,\beta)$ that can be defined recursively as follows.    If $ ind(a) = \mu (\beta)$, then $\cm^+(a, \beta)$
is compact, and so $\cm^+(a, \beta) = \bcm^+(a, \beta)$.   Recursively, suppose  
$\bcm^+(a', \beta')$ for all $a' \in \cp_L$ and $\beta' \in \cp (H)$ with index difference $ind (a') - \mu (\beta') < q$.  Now suppose $a \in \cp_L$ and $\beta \in \cp (H)$ have 
$ind (a) - \mu (\beta) = q$.  Then  $\bcm^+(a, \beta) = \cm^+(a, \beta) \cup \p \bcm^+(a,b)$, where

\begin{equation}\label{bcmplus}
\p \bcm^+(a, \beta)  =    \bigcup_{a \geq a' \in \cp_L}\bigcup_{\beta' \geq \beta \in \cp (H)} \bcm_{L}(a, a') \times \bcm^+ (a', \beta') \times \bcm (\beta', \beta)
\end{equation}
where in these unions, at least one of the inequalities,  $a \geq a'$, and $\beta' \geq \beta$ must be a strict inequality.   Here the space $\bcm_{L}(a, a')$  is the compactification of the space of gradient flows (i.e ``piecewise flows") of the Morse function, $\ce : LM \to \br$ that connect the critical points $a$ and $a'$. The spaces $\bcm (\beta', \beta) = \bcm (\beta', \beta; J, H)$ are as studied above.

We will need another gluing assumption regarding these moduli spaces,  completely analogous to Gluing Assumption \ref{glue}.

\begin{glue}  \label{glue2}     For each $a > b \in \cp_L$ and $\gamma > \beta \in \cp (H)$,  there is   an $\eps >0$ and  local diffeomorphisms 
$$
\bcm_{L}(a, b) \times \bcm^+ (b, \beta)  \times [0, \eps)  \to \bcm^+(a, \beta) \quad \text{and} \quad \bcm^+ (a, \gamma) \times \bcm(\gamma,  \beta)  \times [0, \eps)  \to \bcm^+(a, \beta)
$$ 
which gives $ \bcm^+(a, \beta)$ the structure of a smooth, compact manifold with corners.
\end{glue}

We remark that  the moduli spaces $\bcm^+(a, \beta)$  actually have the structure of $\langle ind(a) - \mu (\beta) \rangle$-manifolds.  The   the faces are given by  
$$
F_j (\cm(a, \beta)) = \coprod_{ind (b) = ind (a)-j}  \bcm_{L}(a, b) \times \bcm^+ (b, \beta)   \quad \sqcup  \coprod_{\mu(\gamma) =   ind (a)-j}\bcm^+ (a, \gamma) \times \bcm(\gamma,  \beta) 
$$

 Notice that the above decomposition of $  \p \bcm^+(a, \beta)$ has two interesting strata.    One is 
\begin{equation}\label{piece}
\coprod_{ind (a) = \mu(\alpha)}\cm^+(a, \alpha) \times \bcm (\alpha, \beta)
\end{equation}
where each $\cm^+(a, \alpha)$ in the above disjoint union is a compact zero dimensional manifold, (and therefore a finite set).  The other stratum of interest is 
$$
\coprod_{ind (b) = \mu(\beta)} \bcm_{L}(a,b) \times \cm^+(b, \beta)
$$
where each $\cm^+(b, \beta)$ is a finite set.  We now prove essentially prove that $\bcm^+(a, \beta)$ defines a framed bordism between these two strata.  More specifically we prove  the following (assuming the two gluing assumptions above).

\begin{lemma}  For each $a \in \cp_L$, $\beta \in \cp (H)$,  the space $\bcm^+(a, \beta)
$ is a compact, framed manifold with corners, in such a way that the framing, when
restricted to the stratum of the boundary given by (\ref{piece}) is given by the product
of framings on $\cm^+(a, \alpha)$ (i.e an orientation on the finite set), and the canonical framing on $\bcm (\alpha, \beta)$  used in  Theorem(\ref{chframed}) above.
\end{lemma}

\begin{proof}  The proof will be very similar to the proof of Theorem \ref{chframed}.  Recall from \cite{abboschwarz} that $\cm^+(a, \beta)$ can be viewed as the space of zeros of a certain section of a vector bundle.  Specifically,  let $\cb^+(a, \alpha)$ be the set of maps $u : \br_+ \times S^1 \to T^*M$ which are Sobolev of type $W^{1,r}$ on every compact subset, and such that
\begin{enumerate}
\item $p\circ u (0, \cdot) \in W^u (a),$  \notag \\
\item there is an $s_0 \geq 0$ for which
$ 
u(s,t) = exp_{x(t)}(\zeta (s,t)) \quad \text{for all} (s,t) \in (s_0, \infty) \times S^1,
$ 
where $\zeta$ is a $W^{1,r}$ section of $\beta^*(TT^*M).$ \notag
\end{enumerate}
There is a vector field on $\cb^+(a, \beta),$ $$\p^+_{J,H} : u \la  \left(\p_su -J(t,u)(\p_tu -X_H(t,u))\right)
$$  whose zeros are $\cm^+(a, \beta)$.  Then, as in the proof of Theorem \ref{chframed}, (see (\ref{stabletan})) the assignment  $u \to D_{fib}\p^+_{J,H}(u)$ defines a homotopy class of map
 
\begin{equation} 
D_{fib}\p^+_{J,H} : \cb^+(a, \beta)  \to BO,
\end{equation} which, when one composes with the inclusion $\cm^+(a, \beta) \hk \cb^+(a, \beta)$
 classifies the stable tangent bundle.  We therefore must show that this map is null homotopic.  Notice that $\cb^+(a,\beta)$ can be viewed as a subspace of $\Omega_{\alpha, \beta}(LT^*M)$, where $\alpha \in \cp (H)$ is the Legendre transform of $a \in \cp_L$.  This can be seen as follows. Let  $u \in \cb^+(a, \beta)$.  Let $u_0 = u(0, \cdot) : S^1 \to T^*M$.
Since $p\circ u_0 : S^1 \to M$ lies in $W^u(a)$, there exists a unique map $v : (-\infty, -1] \to LM$,  which is a gradient  trajectory of $\ce$, and satisfies the initial condition,
$v(-1) = p\circ u_0$.   Let $\tilde v : (-\infty, -] \to T^*M$ be the Legendre transform of $v$.
 That is, it is the composition of $v$ with the map $LM \to LT^*M$ that sends $\gamma$ to $(\gamma, g (\frac{d\gamma}{dt}))$.    Now $\tilde v(-1)$ and $ u_0 : S^1 \to T^*M $ have the same image under composition
 with $p : T^*M \to M$.  So one can define a linear combination,
 $$
 h : [-1, 0] \to L(T^*M)
 $$
 by $ h(t) = -t \tilde v(-1)  + (1+t)u_0.
 $   Now define
$$
 \tilde u : \br   \to L(T^*M)     \quad \text{by} \quad 
 \tilde u(s)  = \begin{cases} \tilde v(s)  \quad \text{for} \quad s < -1 \\
 h(s) \quad \text{for} \quad  -1 \leq s \leq 0 \\
 u(s,\cdot) \quad \text{for} \quad 0 \leq s.
 \end{cases}
 $$

 Then $\tilde u \in \Omega_{\alpha, \beta}(LT^*M)$.   This defines an inclusion map,
 $ 
 \iota : \cb^+(a, \beta) \hk \Omega_{\alpha, \beta}(LT^*M).
 $ 
Furthermore it is clear that $\iota$ is a homotopy equivalence.    Now like  diagram (\ref{diagram}) above, the following diagram homotopy commutes:

\begin{equation}\label{diagram2}
\begin{CD}
\cm^+(a,\beta)   &   @>D_{fib}\p^+_{J,H} >>    &  \bz \times BO \\
@V VV &&  @AAc A \\
\cb^+ (a,\beta)  \xr{\simeq} \Omega_{\alpha, \beta}(L(T^*M))  @>>\Omega L (\tau )>  \Omega L BU(n)  @>>ind \, \bar \p >  \bz \times BU  
\end{CD}
\end{equation}
Here $\tau : T^* M \to BU(n)$ classifies that almost complex tangent bundle of $T^*M$.      
 But by diagram (\ref{tangent}) above, this diagram is homotopic to
 
 \begin{equation}\label{tangent3}
\begin{CD}
\cm^+(a,\beta) @>D_{fib}\p^+_{J,H} >>     \bz \times BO \\
@VVV   @VV \simeq V \\
\Omega L(T^*M) @>>\Omega \rho >  \Omega (U/O)
\end{CD}
\end{equation}

Thus the null homotopy of $\rho : L(T^*M) \to U/O$ defines a null homotopy of 
$D_{fib}\p^+_{J,H} :  \cm^+(a,\beta) \to \bz \times BO$,  and thus a framing.
Notice that since these framings  are induced  applying
the loop functor to a null homotopy of $\rho$,   they are multiplicative, and that in analogy to 
  (\ref{multiply}) above,   they extend to give framings
of the compactified spaces, $\bcm^(a, \beta)$.   Furthermore, when restricted to the 
boundary strata  of the form
$$
\cm^+(a, \alpha) \times \cm(\alpha, \beta) \hk \p (\bcm^+(a, \beta)
$$
these framings are the product of the framing  on $\cm^+(a, \alpha)$, and a framing on $\bcm(\alpha, \beta)$.  The induced framing on $\bcm(\alpha, \beta)$ is given by the null homotopy of the composition
 $$
 \bcm (\alpha, \beta) \hk \Omega_{\alpha, \beta}L(T^*M) \xr{\Omega \rho} \Omega (U/O) \simeq \bz \times BO
 $$
 induced by the canonical null homotopy of $\rho$.  This is   same  framing as the canonical framing used in Theorem \ref{chframed}.   
\end{proof}

\med

Now suppose $a \in \cp_L$ has Morse index $p$, and $\beta \in \cp (H)$ has Conley-Zehnder index $q$.  The above  lemma implies there is a framed embedding
$$
\begin{CD}
\bcm^+(a,\beta) \hk \bcm^+(a,\beta) \times \br^L  @>e > \hk >   J(p ,q-1) \times \br^L
\end{CD}
$$
whose intersection with $J(p,q) \times \br^L \subset  J(p,q-1) \times \br^L$ \, is \, $\coprod_{\mu(\alpha) = p} \cm(a, \alpha) \times \bcm(\alpha, \beta) \times \br^L$.  Similarly, the intersection of the image of $e$  with $J(p-1, q-1) \times \br^L$ \, 
is \, $\coprod_{ind(b) = q} \bcm_L(a, b) \times \cm(b, \beta) \times \br^L$.  This means
that the Pontrjagin-Thom construction,
$$
\zeta_{a, \beta} : J(p, q-1)^+ \wedge S^L_a  \to \left(\bcm^+(a, \beta) \times \br^L\right)\cup \infty \to S^L_\beta
$$
has the property that when restricted to $J(p,q)^+ \wedge S^L_a$
is given by the composition
\begin{equation}\label{comp1}
 \Phi^0_{a, \beta} : J(p,q)^+ \wedge S^L_a \xr{1 \wedge (\bigvee_{\mu(\alpha = p)} \zeta_{a, \alpha})}  J(p,q)^+ \wedge \left(\bigvee_{\mu(\alpha) = p}  S^L_{\alpha}\right)  \xr{\bigvee_\alpha Z_{\cc_H}(\alpha, \beta)} S^L_\beta.
 \end{equation}
 Here  $Z_{\cc_H}(\alpha, \beta)  : J(p,q)^+ \wedge S^L_\alpha \to S^L_{\beta}
 $ is the induced map on morphisms of  the functor $Z_{\cc_H}$ as was defined   (\ref{zalbe}) above, in the proof of Theorem \ref{fltype}.    Notice also that since $ind (a) = \mu (\alpha)$, the map 
 $\zeta_{a, \alpha} : S^L_a \to S^L_\alpha$ is a map whose degree is equal to $\#\cm^+(a, \alpha)$, where the count is defined by the framing of these zero dimensional moduli spaces.  Similarly, when restricted to $J(p-1, q-1)^+ \wedge S^L$, the map $\zeta_{a, \beta}$ is given by the composition

\begin{equation}\label{comp2}
 \Phi^1_{a, \beta} : J(p-1,q-1)^+ \wedge S^L_a \xr{ \bigvee_{ind (b) = q)} Z_{\cc_L}(a, b)}   \bigvee_{ind (b) = q}  S^L_{b}  \xr{\bigvee_b \zeta_{b, \beta}} S^L_\beta.
 \end{equation}
 Here, $\cc_L$ is the flow category of the Morse function $\ce : LM \to \br$ defined in terms of the Lagrangian $L$,  which is a smooth, compact, framed category of Morse-Smale type, where the framings of the morphism spaces $\bcm_L(a,b)$ are induced by the framings of $\bcm^+(a, \beta)$ as above.  This induces a functor $Z_{\cc_L} : \calj \to  Sp_*$ as described in section 1.  As shown in \cite{cjs}, the Floer homotopy type of such a Morse function, given by  the geometric realization, $|Z_{\cc_L}|$ is given by the suspension spectrum of the underlying manifold,
 \begin{equation}\label{morse}
 |Z_{\cc_L}| \simeq \Sigma^\infty (LM_+)   
 \end{equation}
 
 We now show that the  two maps $ \Phi^0_{a, \beta}$ and $ \Phi^1_{a, \beta}
$  are canonically homotopic, via a homotopy that preserves the corner structure. As we will see, these homotopies will patch together to give a homotopy equivalence of the Floer homotopy types given from the Morse function
 $\ce :LM \to \br$, and from the symplectic action, $\ca_H : L(T^*M) \to \br$,
 $$
|Z_{\cc_L}| \simeq |Z_{\cc_H}(T^*M)| = Z(T^*M).
 $$
 Given the equivalence, (\ref{morse}),  we will then know that the Floer homotopy type of $T^*M$ is stably equivalent to the free loop space, as asserted in the theorem. 
  
  To define these homotopies,  let $e_i \in J(p, q-1)$ be the vector with a $1$ in the $i^{th}$ slot and zero's elsewhere,  for each $p > i > q-1$.  Then
 $J(p,q) \subset J(p, q-1)$ is the hyperplane (with corners) given by nonnegative linear combinations of  the $e_i$'s, for  $i = p-1, \cdots q+1 $.   Similarly,  $J(p-1, q-1) \subset J(p, q-1)$is the hyperplane given by nonnegative linear combinations of the $e_i$'s, for
 $i = p-2,  \cdots, q$.    For $t \in [0,1]$, let $J_t \subset J(p, q-1)$ be the hyperplane spanned by nonnegative linear combinations of $(1-t)e_{p-1}+te_q$, and the $e_i$'s for 
 $i = p-2, \cdots , q+1$.    Notice that $J_0 = J(p, q)$, and $J_1 = J(p-1, q-1).$  The restriction of $\zeta_{a, \beta} : J(p, q-1)^+ \wedge S^L_a \to S^L_b$ to $J_t^+ \wedge S^L_a$  is a map we call
\begin{equation}\label{htpy} 
 \Phi^t_{a, \beta} : J_t^+ \wedge S^L_a \la S^L_\beta
\end{equation}
 which gives a canonical homotopy between
 $ \Phi^0_{a, \beta}$ and  $ \Phi^1_{a, \beta}$ viewed as maps 
 $$
  \Phi^0_{a, \beta}, \,  \Phi^1_{a, \beta}: c^{p-q-1}S^L_a \la S^L_{\beta}.
  $$
  Furthermore, this homotopy preserves the corner structure (i.e the boundary stratification of the iterated cones).  
  
  We now define the homotopy equivalence $$\Psi : |Z_{\cc_L}| \xr{\simeq} |Z_{\cc_H}(T^*M)|.$$  It is enough to give a homotopy equivalence, between $|Z^{p,q}_{\cc_L}|$
  and $|Z^{p,q}_{\cc_H}|$ for every $p \geq q$. 
  
  Recall that  for an integer $m$ with $p \geq m \geq q$, $Z_{\cc_L}(m) = \bigvee_{ind (c) = m}S^L_c$.  Similarly, $Z_{\cc_H}(m) = \bigvee_{\mu (\gamma) = m}S^L_\gamma$.
  Define 
  $$
  \Psi_m : Z_{\cc_L}(m) \xr{\simeq} 
Z_{\cc_H}(m) 
$$
to be given by the wedge of the maps $\zeta_{c, \gamma} : S^L_c \to S^L_\gamma$.  Notice that when one applies homology, $\tilde H_*(Z_{\cc_L}(m)) $ is the free abelian group generated by the index $m$  critical points  $\cp_L$ of the Morse function $\ce$.  This is the $m^{th}$ chain group in the Morse complex of the function $\ce$, which we call   $C^{\ce}_*(LM)$.  Similarly,   $\tilde H_*(Z_{\cc_H}(m))$ is the free abelian group generated by the Conley-Zehnder index $m$  critical points $\cp (H)$  of the perturbed symplectic action, $\ca_H$, which is the $m^{th}$ Floer chain group, $CF_m(T^*M)$. Moreover  
\begin{equation}\label{chiso}(\Psi_m)_*: C^{\ce}_m(LM) \to CF_m(T^*M)
\end{equation} is the isomorphism constructed by Abbondondolo and Schwarz in \cite{abboschwarz}.
 
We now consider the induced map
$$
\Psi_m \wedge 1 : Z_{\cc_L}(m) \wedge J(m, q-1)^+ \la  
Z_{\cc_H}(m) \wedge J(m, q-1)^+.
$$
Now recall the relations involved in defining the geometric realization of a functor
$Z : J \to  Sp_*$ (see (\ref{zrealize}).)  For each $m > n > q$ One identifies the image of the embedding
$$
\iota_{m,n}: Z(m) \wedge J(m,n)+ \wedge J(n, q-1)^+ \hk  Z(m) \wedge J(m,q-1)^+
$$
with the image of
$$
j_{m,n} : Z(m)
 \wedge J(m,n)+ \wedge J(n, q-1)^+  \xr{Z(m,n)\wedge 1} Z(n) \wedge J(n, q-1)^+.
 $$
Thus in order for the maps $\Psi_m \wedge 1$ to directly fit together to give a map of geometric realizations, we would need to have that for every $m > n > q$,  $$\Psi_n \circ Z_{\cc_L}(m,n) : Z_{\cc_L}(m) \wedge J(m,n)^+ \to Z_{\cc_L}(n) \to Z_{\cc_H}(n) $$ would need to equal
$$Z_{\cc_H}(n) \circ (\Psi_m \wedge 1): Z_{\cc_L}(m) \wedge J(m,n)^+ \to Z_{\cc_H}(m) \wedge J(m,n)^+ \to Z_{\cc_H}(n).$$  This, however is not true, but these maps are canonically homotopic, via the homotopy $$\Phi^t(m,n) = \bigvee_{\substack{ind (a) = m \\ \mu (\beta) = m}}\Phi^t_{a, \beta}.$$  We accommodate these homotopies by using double mapping cylinders.  Recall that the double mapping cylinder of maps
$
\iota : A \to B \quad \text{and} \quad j : A \to C
$ 
is $$M(\iota, j) = (B \sqcup C) \cup_{\iota, j}A \times I$$ where $(a,0)\in A \times I$ is identified with $\iota (a) \in B$ and $(a, 1) $ is identified with $j(a) \in C$.   There is a projection onto the pushout, which collapses the cylinder,
$$
\pi : M(\iota,j) \to B\cup_A C
$$
which is a homotopy equivalence if $\iota$ is a cofibration.     

Now consider the space $|\tilde Z^{p,q}_{\cc_L}|$ defined like $|Z^{p,q}_{\cc_L}|$
except that rather than identifying the images of $\iota_{m,n}$ and $j_{m,n}$ as  above,
one takes the mapping cylinders, $M(\iota_{m,n}, j_{m,n})$.  By collapsing the cylinders, one gets a map
$$
\pi : |\tilde Z^{p,q}_{\cc_L}| \to |Z^{p,q}_{\cc_L}|
$$
which is a homotopy equivalence, because all of the maps $\iota_{m,n}$ are cofibrations.
 Moreover, the maps $\Psi_m \wedge 1 : Z_{\cc_L}(m) \wedge J(m, q-1)^+ \la  
Z_{\cc_H}(m) \wedge J(m, q-1)^+$ do extend to give a map
$$
\tilde \Psi :  |\tilde Z^{p,q}_{\cc_L}| \to  |Z^{p,q}_{\cc_H}|
$$
by putting the homotopies $\Phi^t(m,n)$ on the cylinders.  Thus we can define
$$\Psi :   |Z^{p,q}_{\cc_L}| \to  |Z^{p,q}_{\cc_H}|$$
to be $\tilde \Psi \circ \pi^{-1}$.  This is well defined up to homotopy.  By (\ref{chiso}) above,  on the level of chain complexes this gives the Abbondandolo-Schwarz isomorphism between the Morse homology and the Floer homology.  Thus
$\Psi$ is a homotopy equivalence.  Moreover by letting $p$ and $q$ vary this defines
a homotopy equivalence,
 $$
 \Psi : \Sigma^\infty (LM_+) \simeq |Z_{\cc_L}| \to |Z_{\cc_H}(T^*M)| = Z(T^*M)
 $$
 as claimed.


\begin{thebibliography}{99}{ 
    \bibitem{abboschwarz}A. Abbondandolo and M. Schwarz, \emph{On the Floer homology of cotangent bundles}, Comm. Pure Appl. Math \bf 59, \rm 254-316 (2006) preprint: math.SG/0408280

\bibitem{abboschwarz2} A. Abbondandolo and M. Schwarz, , \emph{Notes on Floer homology and loop space homology}, 
proceedings of SMS/NATO Adv. Study Instiute on Morse theoretic
methods in non-linear analysis and  symplectic topology, Nato Science Series II: Mathematics, Physics, and Chemistry \textbf{vol 217}, Springer, (2006), 75-108.      

\bibitem{atiyah}M.F.  Atiyah,  \textbf{$K$-theory}, W.A. Benjamin, Inc. (1967).

\bibitem{atiyah2}
  M.F. Atiyah
\emph{Bott Periodicity and the index of elliptic operators},
  Quart. Jour. Math, 
\textbf{19}, (1968), 113 - 140.

\bibitem{baraudcornea}) J.-F Barraud and O. Cornea,  \emph{Lagrangian intersections and the Serre spectral sequence},  preprint (2004),  arXiv: mathDG/0401094

\bibitem{chassull} M. Chas and D. Sullivan, \emph{String topology},  Annals of Math., to appear, arXiv: mathGT/9911159.

\bibitem{cielelai} K. Cieliebak and J.Latschev,  \emph{The role of string topology in symplectic field theory}, preprint (2007), arXiv:mathSG0706.3284

\bibitem{cohenmont} R.L. Cohen, \emph{Morse theory, graphs, and string topology}, 
proceedings of SMS/NATO Adv. Study Instiute on Morse theoretic
methods in non-linear analysis and  symplectic topology, Nato Science Series II: Mathematics, Physics, and Chemistry \textbf{vol 217}, Springer, (2006), 149-184.      

 \bibitem{cg} R.L. Cohen and V. Godin,  \emph{A polarized view of string topology}  Topology, Geometry, and Quantum Field Theory, London Math. Soc.
 Lecture Notes  \bf vol. 308 \rm (2004), 127-154

\bibitem{cj} R.L Cohen and J.D.S. Jones, \emph{A homotopy theoretic realization 
of string topology},  
   \sl Math. Annalen, \bf  vol. 324,  \rm 773-798 (2002).

\bibitem{cjs}R.L. Cohen, J.D.S. Jones, and G.B. Segal, \emph{Floer's infinite dimensional Morse theory and homotopy theory}
  \sl Floer Memorial Volume, \rm
Birkhauser Verlag Prog.
in Math. \bf vol. 133 \rm  (1995), 287 - 325.

 
\bibitem{cjscat} R.L. Cohen, J.D.S. Jones, and G.B. Segal, \emph{Morse theory and classifying spaces}, Stanford University preprint, (1995)  available at http://math.stanford.edu/\~{}ralph/papers.html 

\bibitem{donaldson}S. K. Donaldson,   \bf Floer Homology groups in Yang-Mills theory \rm  with the assistance of M. Furuta and D. Kotschick,  \sl Cambridge Tracts in Mathematics \bf 147, \rm Cambridge Univ. Press, (2002).  

\bibitem{floer} A. Floer,  \emph{Morse theory for Lagrangian intersections} J. Differential Geometry \bf 28 \rm (1988),  513-547.

 \bibitem{fukaya} K. Fukaya, \emph{Application of Floer homology of Lagrangian submanifolds to symplectic topology}, 
proceedings of SMS/NATO Adv. Study Instiute on Morse theoretic
methods in non-linear analysis and  symplectic topology, Nato Science Series II: Mathematics, Physics, and Chemistry \textbf{vol 217}, Springer, (2006), 231-276.      

\bibitem{janich}K. Janich, \emph{On the classification of $O(n)$-manifolds}, Math. nnalen, \bf 176 \rm (1968), 53-76.
 
\bibitem{laures}G. Laures, \emph{On cobordism of manifolds with corners}, Trans. of AMS \bf 352 \rm no. 12,  (2000), 5667-5688.

\bibitem{pressleysegal} 
 A. Pressley and G. Segal,
\textbf{ Loop Groups},
  Oxford Math. Monographs, Clarendon Press
(1986).



\bibitem{salamon} D. Salamon, \emph{Lectures on Floer homology}, Symplectic geometry and topology (Y. Eliashberg and L. Traynor, eds.), IAS/Park City Mathematics Series, AMS (1999), 143-225.

\bibitem{Salzehn}D. Salamon and E. Zehnder, \emph{Morse theory for periodic solutions of Hamiltonian systems and the Maslov index}, Comm. Pure Appl. Math \bf 45 \rm (1992), 1303-1360.

      \bibitem{salamonweber} D. Salamon and J. Weber, \emph{Floer homology and the heat flow}, preprint: math.SG/0304383
      
      \bibitem{segal} G.B. Segal,  \emph{Classifying spaces and spectral sequences}
     Inst. Hautes ƒtudes Sci. Publ. Math. No. 34 1968 105--112. 
     
     \bibitem{taubes}  C. Taubes, $L^2$\bf - moduli spaces on $4$-manifolds with cylindrical  ends,  \rm \sl Monographs in Geometry and Topology, I, \rm International Press (1993).

  \bibitem{viterbo}C. Viterbo, \emph{Functors and computations in Floer homology with applications, Part II}, preprint, (1996). 
    
     \bibitem{weber}J. Weber, \emph{Perturbed closed geodesics are periodic orbits: Index and transversality}, Math. Zeit. {\bf (241)}, (2002), 45-81.

     }\end{thebibliography}
 \end{document}